\documentclass{article}
\usepackage{amsmath}
\usepackage{amsthm,verbatim}
\usepackage[all]{xy}

\newtheorem{lmm}{Lemma}[section]

\newtheorem{prp}{Proposition}[section]

\newtheorem{thm}{Theorem}[section]

\newtheorem{dfn}{Definition}[section]
\newtheorem{rem}{Remark}%%%

\newcommand{\Z}{{\mathbf Z}}                   %Entiers
\newcommand{\R}{{\mathbf R}}                   %Réels
\newcommand{\C}{{\mathbf C}}                   %Droite complexe
\newcommand{\Ct}{\widehat{\mathbf C}}              %Droite complexe duale
\newcommand{\CP}{\mathbf{CP}^1}                %Sphere de Riemann
\newcommand{\CPt}{\widehat{\mathbf{CP}}^1}     %Sphere de Riemann duale
\newcommand{\D}{\mathcal{D}}                   %Opérateurs différentiels holomorphes
\newcommand{\DM}{\mathcal{M}}                  %D-module du fibré méromorphe
\newcommand{\DN}{\mathcal{N}}                  %D-module général
\newcommand{\DMt}{\widehat{\mathcal{M}}}           %D-module transformé minimal
\newcommand{\Comp}{\mathcal{C}}                %Complexe simple associé à (\ref{compldouble})
\newcommand{\Lisse}{\mathcal{C}^{\infty}}           %C^infini
\newcommand{\dbar}{\bar{\partial}}                    %d-bar operator
\newcommand{\ra}{\to}                  %fleche
\newcommand{\lra}{\longrightarrow}             %fleche longue
\newcommand{\Om}{{\mathbf \Omega}}                  %formes holomorphes
\newcommand{\Hol}{\mathcal{O}}                 %fonctions holomorphes
\renewcommand{\d}{\mbox{d}}                    %derivation exterieure
\newcommand{\dt}{\widehat{\mbox{d}}}               %derivation triviale dans la direction de xi
\newcommand{\Et}{\widehat{E}}                      %fibre transforme
\newcommand{\Nt}{\widehat{\nabla}}                      %connexion transformee
\newcommand{\Pt}{\widehat{P}}                      %lieu singulier transforme
\newcommand{\trmetr}{\widehat{h}}                      %metrique harmonique transforme
\newcommand{\ti}{\widehat{\infty}}                    %point a l'infini transforme
\newcommand{\Mt}{\widehat{M}}                      %module transforme
\newcommand{\M}{{\mathbf M}}                   %module en deux variables
\renewcommand{\H}{{\mathbf H}}                   %hypercohomologie
\newcommand{\A}{\mathcal{A}}                 %famille des valeurs propres du D-module
\newcommand{\isom}{\simeq}                    %isomorphisme
\newcommand{\gr}{\mbox{gr}}                   %gradués
\DeclareMathOperator{\pdeg}{par-deg}                 %degré
\DeclareMathOperator{\parab}{par-}
\DeclareMathOperator{\rank}{rk}
\DeclareMathOperator{\res}{Res}
\DeclareMathOperator{\im}{Im}
\DeclareMathOperator{\coker}{coker}
\DeclareMathOperator{\diag}{diag}
\DeclareMathOperator{\Gl}{Gl}

\title{Nahm transform and parabolic minimal Laplace transform}

\author{Szil\'ard \textsc{Szab\'o}
          \footnote{Department of Geometry, Mathematical Institute, Faculty of Science,
            Budapest University of Technology and Economics,
            Egry J. u. 1, H \'ep., H-1111 Budapest, Hungary, {\tt szabosz@math.bme.hu}
}}
\begin{document}

\maketitle

\tableofcontents

\begin{abstract}
We prove that Nahm transform for integrable connections with a finite number of regular 
singularities and an irregular singularity of rank 1 on the Riemann sphere is equivalent 
-- up to considering integrable connections as holonomic $\D$-modules -- to minimal Laplace transform. 
We assume semi-simplicity and resonance-freeness conditions, and we work in the framework 
of objects with a parabolic structure. In particular, we describe the definition of the 
parabolic version of Laplace transform due to C. Sabbah. The proof of the main result relies 
on the study of a twisted de Rham complex.
\end{abstract}

\section*{Introduction}

Nahm transform is a correspondence in the theory of $4$-manifolds, between solutions of the 
anti-self-dual (ASD) Yang-Mills equations on $\R^{4}$, invariant by a closed additive subgroup 
and with finite energy on the quotient on the one hand, and solutions of the ASD equations on 
the dual vector space  $(\R^{4})^{\ast}$, invariant with respect to the dual subgroup and 
with finite energy on the quotient on the other hand. It can be considered as a differential 
geometric variant of Fourier-Mukai transform \cite{bbh}. 

%$Various cases of this transformation have been studied by Atiyah-Drinfeld-Hitchin-Manin, 
%Nahm, Mukai, Schenck, Braam-van Baal, Cherkis-Kapustin, Jardim, 
%Biquard-Jardim, Bonsdorff, Nye, Hurtubise-Charbonneau, etc. 
When the first subgroup is $\R^{2}$, and after identification of the quotient $\R^{2}$ with 
$\C$, invariant solutions of the ASD equations are holomorphic connections together with a 
harmonic metric (see \cite{hit}). 
On the other hand, the dual quotient is then again a copy of $\C$, that we shall denote by 
$\Ct$ in this text to avoid confusion. Hence, Nahm transform maps holomorphic connections 
with a harmonic metric on $\C$ into holomorphic connections 
with a harmonic metric on $\Ct$. 
In \cite{Sz}, we describe this transform for holomorphic connections on the Riemann sphere 
$\CP$ with a finite number of regular singularities in $\C$, and an irregular singularity 
of Poincar\'e rank $1$ at infinity, and endowed with a parabolic structure in all singular points. 
It is well-known that holomorphic connections can be viewed as a special case of holonomic 
$\D$-modules of finite rank. 
%and under some assumptions this correspondence is an equivalence of categories (see Proposition \ref{prop:equiv}). 
Furthermore, this correspondence is compatible with parabolic structures. 
On the other hand, there exists a so-called minimal (Fourier-)Laplace transform for $\D$-modules 
without a parabolic structure (see e.g. \cite{malgr}). 
Claude Sabbah extends in \cite{sab} this construction to the case with parabolic structure 
in the singularities. 
We shall review this extension in Section \ref{sec:parlapl}. 
He also asked whether Nahm transform and parabolic minimal Laplace transform agree under the 
above-mentioned equivalence of categories. Our main result is an affirmative answer to this 
question (see Theorem \ref{thm:main}). This result allows us to give a complex algebraic 
definition of Nahm transform, which is initially defined using $L^2$-theory. An immediate 
consequence is that Nahm transform is a holomorphic map between moduli spaces of 
stable holomorphic connections. This is one step in showing that Nahm transform is actually a 
hyper-K\"ahler isomorphism between moduli spaces. 

A few words are due here to explain our assumptions on the singularities 
(regular singularities except for one irregular singularity of Poincar\'e rank $1$ at infinity). 
At first sight this might look an artificial choice; however, it turns out 
that this class of connections is preserved by our transformation, 
hence it is natural to study the properties of this correpondence. 
More importantly, this class of connections appears naturally in the theory of 
Frobenius manifolds \cite{sab2} as well as in Mirror Symmetry \cite{kkp}. 
Notice however that Laplace transform is defined not only under these 
assumptions; indeed, in a forthcoming paper \cite{biq-sz} we plan to extend 
the definition of Nahm transform for 
integrable connections with harmonic metric on the Riemann sphere to 
a more general setup, and we hope that the present results can be 
extended to that situation too.

\paragraph*{Acknowledgments} The author would like to express his thanks to Claude Sabbah for 
drawing his attention on this question, and for the useful remarks during the preparation 
of this text. 
Various parts of this text were written while the author stayed at several
institutes. It was started during his scholarship at the Max Planck Institute for Mathematics 
in Bonn: he would like to thank for the hospitality and excellent working conditions.
Later parts were written while the author was a junior researcher fellow at the 
Alfr\'ed R\'enyi Institute of Mathematics in Budapest; he wishes to thank for the grant provided. 
Thanks equally to the Budapest University of Technology and Economics, where the text was finished. 

\section{Parabolic connections and Nahm transform}\label{sec:nahm}
We briefly explain the first correspondence. 
Let $\C$ denote the complex line with its canonical holomorphic coordinate $x=x_{1}+ix_{2}$ 
(denoted $z$ in \cite{Sz}), 
and  $\CP$ stand for its one-point compactification, the Riemann sphere. The point at 
infinity will be denoted $\infty$. 
Let $E$ be a holomorphic vector bundle of rank $r$ on $\CP$, whose underlying $\Lisse$ 
vector bundle and holomorphic structure will be denoted by $\Lisse(E)$ and $\dbar^{E}$ 
respectively. 
Let $P=\{ p_{1}, \ldots , p_{n} \}$ be a fixed finite set in $\C$, the singular locus at 
finite distance. 
Let $\Hol$ and $\Om^1$ stand for the sheaf of holomorphic functions and holomorphic $1$-forms
on $\CP$ respectively. 
We will denote by $\Om^{1}(\ast P \cup \{ \infty \})$ the sheaf of meromorphic $1$-forms 
with poles of arbitrary order in the points of $P$ and at infinity, and no other poles. 
Let $\nabla$ be a meromorphic connection on $E$ with poles in $P\cup \{ \infty \}$: 
\begin{equation*}
      \nabla : E \lra \Om^{1} (\ast P \cup \{ \infty \}) \otimes_{\Hol} E.
\end{equation*}
We suppose that $\nabla$ is  logarithmic in the points of $P$, i.e. in a local holomorphic 
trivialization of $E$ near a point $p_j\in P$, it can be written 
\begin{equation*}
     \d + M^{j} \d x 
\end{equation*}
with $M^{j}$ a matrix whose coefficients have at most a simple pole.  
We also suppose that the residue of $\nabla$ in all $p_j$ is semi-simple, 
and hence in a suitable trivialization it can be written 
\begin{equation*}
    \d + A^{j} \frac{\d x}{x-p_{j}} + \mbox{holomorphic terms},
\end{equation*}
where $A^j$ is a diagonal matrix
\begin{equation}\label{aj}
    A^{j}=\begin{pmatrix}
                  0 & & & & & \\
                  & \ddots & & & & \\
                  & & 0 & & & \\
                  & & & \mu^{j}_{r_{j}+1} & & \\
                  & & & & \ddots & \\ 
                  & & & & & \mu^{j}_{r}
           \end{pmatrix},
\end{equation}
with $\mu^{j}_{k}\neq 0$ for $r_{j}<k\leq r$.
Let $E$ be endowed with a parabolic structure in $P$: this means that there exists a 
decreasing exhaustive filtration $E_{p_{j}}^{\bullet}$ of $E_{p_{j}}$ indexed by a 
finite set of $[0,1[$, compatible with the residue in the sense that all 
$E_{p_{j}}^{\beta}$ is generated by the last $k$ elements of the above diagonalizing 
trivialization of $A^{j}$, for a suitable $k=k(\beta)$ between $0$ and $r$. 
The numbers $\beta^j_k \in [0,1[$ such that the graded pieces 
$\gr^{\beta^j_k}E_{p_{j}}^{\bullet}$ 
of the filtration are non-trivial are called parabolic weights. 

\begin{rem}\label{metrrem}
The parabolic structure usually comes from a harmonic Hermitian metric \cite{sim}. 
Namely, assume given a tame harmonic metric $h$; 
then the filtration is characterised by the requirement that for any $\beta$ 
the vectors $e\in E_{p_{j}}^{\beta}\setminus 0$ are exactly those vectors 
whose holomorphic extensions $\sigma$ satisfy 
%for all $\varepsilon>0$ the norm of $\sigma$ with respect to a harmonic metric $h$ should satisfy the estimate
$$
     h(\sigma,\sigma)  < c|x-p_{j}|^{2\beta-\varepsilon}
$$
for all $\varepsilon>0$ and for a suitable $c>0$ (depending on $e$ and $\varepsilon$).
Under our assumption of semi-simplicity, for any $e\in\gr^{\beta}E_{p_{j}}^{\bullet}\setminus 0$ 
and holomorphic extension $\sigma$ the metric actually satisfies the stronger asymptotic equality 
$$
      h(\sigma,\sigma) \approx  |e|^2|x-p_{j}|^{2\beta}
$$
for some norm $|.|$ on $E_{p_{j}}$ coming from a scalar product, where $\approx$ means that the quotient of the
two sides converges to $1$ as $x \ra p_{j}$. 
In fact, one can define the notion of parabolic structure for connections 
with arbitrary logarithmic singularities at finite distance (not necessarily semi-simple ones),
and Nahm transform probably carries through to this setup too. 
In this case, we also need to consider the weight filtration of the nilpotent part of the 
residue on $E_{p_{j}}$, and the harmonic metrics we are interested in have a more refined 
behavior 
$$
    h(\sigma,\sigma) \approx |e|^2 |x-p_{j}|^{2\beta} \log(|x-p_{j}|)^{w}, 
$$
where $w$ is the weight of $e$ in the weight filtration. 
However, in this paper we will stick to the semi-simple case. 
\end{rem}

The behavior of the singularity at infinity is supposed to be of Poincar\'e rank $1$, 
with diagonalizable polar part: in a suitable trivialization of $E$ near $\infty$, 
it can be written 
\begin{equation*}
      \d + M^{\infty} \d x, 
\end{equation*}
where $M^{\infty}$ is a holomorphic matrix, and such that up to lower-order terms 
\begin{equation}\label{connirr}
   M^{\infty} = A + \frac{C}{x}. 
\end{equation}
Here 
$$
     A =  \begin{pmatrix}
                \xi_{1} & & & & & & \\
                & \ddots & & & & & \\
                & & \xi_{1}  & & & & \\
                & & & \ddots & & & \\
                & & & & \xi_{n'} & & \\
                & & & & & \ddots & \\
                & & & & & & \xi_{n'}
                        \end{pmatrix}
$$ 
and 
$$
    C = \begin{pmatrix}
                \mu_1^{\infty} & & \\
                & \ddots & \\
                & & \mu_r^{\infty}
                        \end{pmatrix}
$$
are constant diagonal matrices, the eigenvalue $\xi_l$ of $A$ being displayed on the 
successive places $1+a_{l},\ldots,a_{1+l}$ of the diagonal. 
Finally, we suppose given a parabolic structure in this singularity too: 
this is again a decreasing exhaustive filtration $E_{\infty}^{\bullet}$ of $E_{\infty}$ 
indexed by a finite set of $[0,1[$ (the parabolic weights $\beta^{\infty}_k$ at infinity), 
compatible with the polar part of the connection. 
Again, the parabolic structure is related to the behaviour of a suitable harmonic metric $h$ 
near infinity: a holomorphic section $\sigma$ extending a non-zero vector in 
$\gr^{\beta^{\infty}_k}E_{\infty}^{\bullet}$ is required to have the asymptotic 
\begin{equation}\label{asyminf}
   h(\sigma,\sigma) \approx |x|^{-2\beta^{\infty}_k}
\end{equation}
up to a constant multiple. 

\begin{rem}
The behavior (\ref{aj}) means in the terminology of Simpson \cite{sim} that $(E,D)$ 
is a regular filtered flat bundle with semi-simple residue near the singular points at 
finite distance. The behaviour at infinity is irregular of quite simple type;
such a connection (with only $0$ as regular singular point, but together with some extra data) 
is called non-commutative Hodge structure of exponential type by Katzarkov, Kontsevich and Pantev 
\cite{kkp}, where one can also find an account of other related definitions existing in the literature. 
\end{rem}

The following notions will be crucial for several statements of our paper. 

\begin{dfn}\label{def:rf}
We say that the connection is \emph{resonance-free} if the following conditions hold: 
\begin{enumerate}
\item{for all $p_{j} \in P$ and $r_{j}<k,m\leq r$ one has $\Re(\mu^{j}_{k})\notin\Z$, 
    $\mu^{j}_{k}\neq \beta^{j}_{k}$, furthermore $\mu^{j}_{k}-\mu^{j}_{m}\notin\Z\setminus\{ 0\}$;}
\item{for all $1 \leq k \leq r$, one has $\Re(\mu^{\infty}_{k}) \notin \Z$, 
    $\mu^{\infty}_{k}\neq \beta^{\infty}_{k}$, furthermore 
    $\mu^{\infty}_{k}-\mu^{\infty}_{m}\notin\Z\setminus\{ 0\}$ for any $a_l < k,m \leq a_{1+l}$.}
\end{enumerate}
\end{dfn}

\begin{rem}
Clearly, these assumptions are fulfilled for generic choices of singularity parameters. 
The first and third of each set of conditions simply express that eigenvalues of the residue do not differ 
by non-zero integers (recall that if $r_j>0$ then $\mu^j_1=\cdots =\mu^j_{r_j}=0$ is also an eigenvalue); 
we included them in order to simplify the treatment of parabolic structures in the 
$\mathcal{D}$-module setup, but we expect that they are not necessary for the results to hold. 
The second one states that the parabolic weights of the associated Higgs bundle are non-zero on the 
image of the residue of the Higgs field (c.f. the Table on page 720 of \cite{sim}); 
we included it because we will rely on results of \cite{Sz} where this is assumed. 
\end{rem}

We will also impose a compatibility condition between the polar part of the connection 
and the parabolic structure. 

\begin{dfn}\label{admhol}
We say that the parabolic structure is \emph{admissible} if $\gr^{0}E_{p_{j}}^{\bullet}$ 
contains the eigenspace $\psi^0_{p_{j}}(E,\nabla )$ of the residue for the 
$0$-eigenvalue; in different terms, the natural map 
$$
    \psi^0_{p_{j}}(E,\nabla )\to \gr^{0}E_{p_{j}}^{\bullet}
$$
is injective (including the possibility $\psi^0=0$). 
\end{dfn}

\begin{rem}
This condition will be needed for the $L^2$-resolution Lemma \ref{lem2.1}. 
Notice that it is weaker than its counterpart in \cite{Sz} in that it also allows zero 
parabolic weight correspond to non-zero eigenvalues of the residue. 
In a forthcoming paper \cite{biq-sz} we plan on extending Nahm transform to a more general setup; 
for example, we will not require there that the residue of the Higgs field be semi-simple and that 
the parabolic weights corresponding to the non-zero eigenvalues of the Higgs bundle be all positive. 
The hope is that the results of this paper will carry through to the more general situation as well. 
\end{rem}

As usual, we call parabolic degree of $(E,\nabla)$ the real number 
\begin{equation} \label{pardeg}
    \pdeg(E)=\deg(E) + \sum_{j\in P \cup \{ \infty \}, k\in \{1,\ldots,r\}} \beta_{k}^{j}, 
\end{equation}
and \emph{parabolic slope} the number 
\begin{equation} \label{parpente}
    \parab\mu(E)=\frac{\pdeg(E)}{\rank(E)}. 
\end{equation}
Moreover, we say that $(E,\nabla)$ is (parabolically) stable if any $\nabla$-invariant 
holomorphic sub-bundle $F$, endowed with the induced parabolic structure, satisfies  
\begin{equation} \label{stable}
     \parab\mu(F) < \parab\mu(E).
\end{equation}
The following result is important for our arguments: 
\begin{thm}[C. Simpson \cite{sim}, C.Sabbah \cite{sab-hi}, O.Biquard-Ph.Boalch \cite{biqboa}]
If $(E,\nabla)$ is a parabolically stable holomorphic connection 
of zero parabolic degree on a punctured compact
curve with arbitrary diagonalizable polar parts, then there exists a unique harmonic 
metric (up to multiplication by a constant) behaving as required in Remark \ref{metrrem} 
and in (\ref{asyminf}). 
\end{thm}
A harmonic metric is one that satisfies a certain second-order non-linear partial 
differential equation that we do not spell out here, see e.g. \cite{hit} or
\cite{sim}. Roughly speaking, it is the differential equation governing 
the fact that the metric, seen as a map from the universal cover of the 
punctured Riemann surface to the symmetric space of Hermitian metrics 
on a fixed vector space, be harmonic as a map between Riemannian manifolds. 

From now on, we will suppose that $(E,\nabla)$ is a stable, resonance-free meromorphic 
connection of zero parabolic degree on $\CP$ having regular singularities with
semi-simple residues in $P$, an irregular 
singularity of Poincar\'e rank $1$ at infinity with diagonalizable polar part, and endowed 
with an admissible parabolic structure in all singularities (at infinity, this 
simply means a parabolic structure without any further assumption). 

Denote by $\Ct$ the dual complex line of $\C$, and by $\CPt$ the one-point compactification 
of $\Ct$.
Given $(E,\nabla)$ as above, we define in \cite{Sz} a holomorphic bundle $\Et$ on $\CPt$ and 
a meromorphic connection $\Nt$ with semi-simple regular singularities in the points 
$\Pt= \{ \xi_{1}, \ldots ,\xi_{n'} \}$, and an irregular singularity of Poincar\'e rank $1$ 
at infinity with diagonalizable polar part, called the Nahm transform of $(E,\nabla)$, with 
an admissible parabolic structure in all singularities. 
In what follows, we outline its construction and main properties. 
First, consider the differential-geometric flat connection defined on $\C\setminus P$ 
\begin{equation*}
      D=\dbar^{E}+\nabla : \Lisse(E) \lra \Omega^{1} \otimes_{\Lisse} \Lisse(E),
\end{equation*}
where $\Omega^{p}$ stands for smooth $p$-forms. Let now $\xi \in \Ct - \Pt$ be a parameter 
and twist $D$ by the formula 
\begin{equation*}
     D_{\xi}=D-\xi \d x. 
\end{equation*}
Consider the elliptic differential complex
\begin{equation}\label{ellcompl}
        \Omega^{0}\otimes E \xrightarrow{\text{$D_{\xi}$}}  \Omega^1 \otimes E 
                  \xrightarrow{\text{$D_{\xi}$}} \Omega^2 \otimes E .
\end{equation}
Endow $\C$ with the standard Euclidean metric; it is singular at infinity. 
It is then possible to show (Theorems 2.6 and 2.16 of \cite{Sz}) that 
the $L^2$-cohomologies of degrees $0$ and $2$ of this complex vanish for 
all $\xi$, and that its first $L^2$-cohomology $L^{2}H^{1}(D_{\xi})$ is a finite-dimensional 
vector space, whose dimension is independent of $\xi$. We let the fiber of $\Et$ 
over the point $\xi$ be $L^{2}H^{1}(D_{\xi})$. As $\xi$ varies, by the implicit function
theorem in Hilbert space, these finite-dimensional subspaces of $L^2(\Omega^1 \otimes E)$ 
define a smooth vector bundle $\Et|_{\Ct - \Pt}$ over $\Ct - \Pt$. 
It inherits a holomorphic structure from the trivial 
holomorphic structure $\dt^{0,1}$ on the trivial Hilbert bundle $L^2(\Omega^1 \otimes E)$, 
and a holomorphic connection $\Nt$ from $\dt^{1,0}- x \d \xi$, where
$\dt^{1,0}$ stands for the $(1,0)$-part of the trivial connection. The fiber $\Et_{\xi}$ 
also has an interpretation as the $L^2$-kernel of the Laplace operator 
\begin{equation}\label{Laplacian}
    \Delta_{\xi} = D_{\xi}D_{\xi}^* + D_{\xi}^* D_{\xi}, 
\end{equation}
where $D_{\xi}^*$ is the adjoint of the connection operator with respect to the harmonic 
metric; in other words as the space of $L^2$ harmonic $1$-forms (Theorem 2.21 of \cite{Sz}). 
The $L^2$-norm of such a $1$-form defines a Hermitian metric $\trmetr_{\xi}$ on the fiber $\Et_{\xi}$, 
and these fiber metrics vary smoothly with $\xi$ over $\Ct - \Pt$. 
This metric is called the transformed metric, and one can prove that if the original 
parabolic structure is admissible, then $\trmetr$ induces an extension of 
$\Et$ as a holomorphic bundle to $\CPt$ defined as the sheaf of local 
holomorphic sections of bounded norm (Section 4.4 of \cite{Sz}). 
The behaviour of $\Nt$ near the singularities with respect to this extension is of the same type as that of 
$\nabla$ described in (\ref{aj}) and (\ref{connirr}) (with different parameters). 
In addition, this extension of $\Et$ also induces an admissible parabolic structure on 
$\Et$ at the points of $\Pt \cup \ti$, i.e. the transformed metric is adapted to the transformed 
parabolic structure (Section 4.6 of \cite{Sz}). 
Furthermore, one can prove that if the original metric $h$ is harmonic, 
then the same thing holds for $\trmetr$ (Theorem 4.9 of \cite{Sz}). 
In particular, as the existence of an adapted harmonic metric is equivalent to poly-stability \cite{biqboa}, 
it follows that if the original connection $(E,\nabla)$ is poly-stable, then so is $(\Et,\Nt)$. 
Finally, we obtain an explicit description of the singularity behavior of $\Nt$ in 
the singular set; in particular, we can show that if $\nabla$ is resonance-free, then so 
is $\Nt$ (Theorems 4.30 and 4.31 of \cite{Sz}). 

Let us now recall in more details how the singularity behavior of $\Nt$ could be obtained. 
We will only work with the singularity at infinity, because the case of 
the finitely located singular points is similar (and in fact, even simpler). 
The computation uses Higgs bundles instead of integrable connections, and 
relies on the formulae of \cite{sim} for the link between the singularity 
data of the two. Hence, it is sufficient to show how the singularity data 
of the corresponding Higgs bundle transforms. For this purpose, 
in Section 4.6 of \cite{Sz} we constructed explicitly a family of 
\emph{asymptotically harmonic} $1$-forms $\hat{\sigma}(\xi)$ in $\xi$, i.e. such that the 
$L^{2}$-norm of $\Delta_{\xi}\hat{\sigma}(\xi)$ is negligible compared to the $L^{2}$-norm of 
$\hat{\sigma}(\xi)$ (the quotient 
$\Vert\Delta_{\xi}\hat{\sigma}(\xi)\Vert_{L^2}/\Vert\hat{\sigma}(\xi)\Vert_{L^2}$ 
converges to $0$ as $\xi \ra \xi_{l}$), and such that 
in a neighborhood of $\ti$ they match up to give a holomorphic local section for the 
holomorphic structure of the Higgs bundle associated to $\widehat{E}$. 
We give an outline of that construction. Call $(\mathcal{E},\theta )$ the Higgs bundle 
associated to $(E,\nabla)$, and set 
\begin{equation}\label{D''}
   D_{\xi}''=\dbar^{\mathcal{E}}+(\theta- \xi \d x/2).
\end{equation}
We call the operator 
$$
  (D'')^j = \bar{\partial}^{\mathcal{E}} + \frac{A^j-\diag(\beta^j_k)_{k=1}^r}{2}\frac{\d x}{x-p_j}
$$
the model $D''$-operator, and we deform it as 
$$
   (D_{\xi}'')^j= (D'')^j -\frac{\xi}{2}\d x. 
$$
We introduce the model Hermitian metric 
$$
   h^j=\diag(|x-p_j|^{2\alpha^j_k})_{k=1}^r, 
$$
with 
$$
   \alpha^j_k = \{ \Re (\mu^j_k) \}, 
$$
where $\{ .\}$ stands for the fractional part in $[0,1[$ of the number. 
Finally, denote by $(D_{\xi}'')^{j*}$ the formal adjoint of $(D_{\xi}'')^j$ for the model 
metric $h^j$, and set 
$$
   (\Delta_{\xi}'')^j = (D_{\xi}'')^{j*}(D_{\xi}'')^j+(D_{\xi}'')^j(D_{\xi}'')^{j*}
$$
for the model Laplace operator of $(D_{\xi}'')^j$. 
Now, define the spectral points corresponding to $\xi$ 
as the set of all $x \in \C$ such that $\det(\theta(x) - \xi \d x/2)=0$. 
Of course, these points vary analytically with $\xi$, and as $\xi \ra \ti$,
they converge to some $p_j$ asymptotically to constant multiples of $\xi^{-1}$. 
Denote by $\chi(r)$ the standard Gaussian function on the plane as a function of the radius, 
cut-off smoothly outside the disk of radius $r>1$ so that it is supported on the disk of radius $2$. 
Let $\sigma^{j}_{k}$ be a holomorphic section of $\mathcal{E}(p_j)$ near $p_j$ 
extending an eigenvector in the image of the residue of the Higgs field, 
such that for each $\xi$ the vector $\sigma^{j}_{k}(q_k(\xi))$ represent a 
non-vanishing class in $\coker(\theta(q_k(\xi)) - \xi \d x/2)$, 
and choose a sufficiently small $\varepsilon_{0}>0$. 
Define the ${E}$-valued smooth $(1,0)$-form on $\CP$ for every $\xi$ close to $\ti$
by the formula 
\begin{equation}\label{trsection}
v^{j}_{k}(\xi,x)\d x=\chi(\varepsilon_{0}^{-1}|\xi|(x-q_{k}(\xi)))\sigma^{j}_{k}(x)\d x;
\end{equation} 
it is supported on a disk of radius $2\varepsilon_{0}|\xi|^{-1}$ centered at the spectral point $q_{k}(\xi)$. 
Then there exists an  ${E}$-valued $(0,1)$-form $t^{j}_{k}(\xi,x)\d \bar{x}$, necessarily supported in 
the same disk, such that 
$$
   (D_{\xi}'')^j(v^{j}_{k}(\xi,x)\d x+t^{j}_{k}(\xi,x)\d \bar{x})=0
$$ 
for all $\xi$. 
We define $\hat{\sigma}^{\infty}_{k}(\xi,x)$ as the harmonic representative of this cohomology class 
with respect to the model Laplace operator, cut off by $\chi$ in a neighborhood of $q_{k}(\xi)$ 
of diameter $O(|\xi|^{-1})$. 
It can then be shown (see Lemma 4.34 of \cite{Sz}) that the family in $\xi$ of 
the sections $\hat{\sigma}^{\infty}_{k}(\xi,z)$ is asymptotically harmonic. 
On the other hand, the $L^{2}$-norm on $\C$ of $\hat{\sigma}^{\infty}_{k}(\xi,z)$ 
is approximately equal to $|\xi|^{-\alpha^{j}_{k}}$ 
multiplied by the $L^2$-norm of a Gaussian function normalized polynomially 
and of at most polynomial height with respect to $|\xi|$, and is therefore bounded 
from above by a polynomial in $|\xi|$. 

\section{Laplace transform without parabolic structure}\label{sec:lapl}
In this section, we recall the classical notions of a module over the Weyl-algebra and 
its Laplace-transform. 

Let ${\C}[x]\langle\partial_{x}\rangle$ denote the Weyl algebra in the variable $x$: 
this is the algebra of polynomial differential operators over ${\C}$, generated by 
the formal variables $x$ and $\partial_{x}$ and the relation $[x,\partial_{x}]=-1$.
Let $M$ be a left ${\C}[x]\langle\partial_{x}\rangle$-module; this means that we 
have an action of $x$ and $\partial_{x}$ on the elements of $M$, satisfying 
$$
   x \cdot (\partial_{x}\cdot m)- \partial_{x}\cdot (x\cdot m)=-m.
$$
In what follows, when we speak of module, we will always mean left-module. 
Let $\xi$ be another variable. Then there exists an isomorphism between the Weyl algebras 
${\C}[x]\langle\partial_{x}\rangle$ and ${\C}[\xi]\langle\partial_{\xi}\rangle$: 
\begin{align}
     {\C}[x]\langle\partial_{x}\rangle & \lra {\C}[\xi]\langle\partial_{\xi}\rangle \notag \\
     x & \mapsto - \partial_{\xi} \label{lapl1} \\
     \partial_{x} & \mapsto \xi \label{lapl2}.
\end{align}
A ${\C}[x]\langle\partial_{x}\rangle$-module $M$ automatically becomes a 
${\C}[\xi]\langle\partial_{\xi}\rangle$-module via this isomorphism; we denote this 
resulting module $\Mt$, and call it the Laplace-transform of $M$. 
In the same way, a ${\C}[\xi]\langle\partial_{\xi}\rangle$-module $N$ becomes 
automatically a ${\C}[x]\langle\partial_{x}\rangle$-module, denoted $\check{N}$ and 
called the inverse Laplace-transform of $N$. These operations are clearly inverses of each other. 

\subsection{Interpretation as a cokernel} \label{conoyauinterpr}
It is well-known that Laplace transform has an interpretation as a cokernel 
(see \cite{kl} Lemma 7.1.4). We recall this interpretation here. Denote by 
\begin{equation*}
    {\C}[x,\xi]\langle\partial_{x},\partial_{\xi}\rangle
\end{equation*}
the Weyl-algebra in the variables $x$ and $\xi$, i.e. the algebra over $\C$ 
generated by the variables $x,\xi,\partial_{x},\partial_{\xi}$ and relations
\begin{align*}
     [x,\partial_{x}]&=-1 \\
     [\xi ,\partial_{\xi }]&=-1,
\end{align*}
all other generators commuting. 
Let the variables $\xi,\partial_{\xi}$ act trivially on $M$; this induces 
a ${\C}[x,\xi]\langle\partial_{x},\partial_{\xi}\rangle$-module 
\begin{equation*}
     \M = M \otimes_{\C} \C[\xi]
\end{equation*}
with $x,\partial_{x}$ acting on $M$ and $\xi,\partial_{\xi}$ acting on $\C[\xi]$ in the natural way. 
Then the kernel of the map 
\begin{equation*}
     \M \xrightarrow{\partial_{x}-\xi} \M
\end{equation*}
vanishes, and its cokernel is bijective with $M$ as a set. 
Moreover, the cokernel inherits the structure of a 
${\C}[\xi]\langle\partial_{\xi}\rangle$-module from $\M$: the action of $\partial_{\xi}$ 
is induced by the operator $\partial_{x}+\partial_{\xi }-x-\xi$ on $\M$. 
Notice that the action of $\partial_{\xi}$ on the cokernel is then 
induced by the action of $\partial_{\xi }-x$ on $\M$ too, for $\partial_{x}-\xi$ clearly 
acts trivially on its own cokernel. Therefore, for any $m \in M$, $\partial_{\xi}$ 
applied to the class of the element $m \otimes 1$ in $\coker(\partial_{x}-\xi)$ 
is $(\partial_{\xi }-x)\cdot (m \otimes 1)=-x\cdot m$, so the action of $\partial_{\xi}$ 
on the cokernel is equal to the action of $-x$ on $M$. Similarly, the class of 
$m \otimes 1$ in $\coker(\partial_{x}-\xi)$ is annihilated by $\partial_{x}-\xi$, so 
the action of $\xi$ on the cokernel is by $\partial_{x}$. Hence, we get the transformation 
given by (\ref{lapl1}-\ref{lapl2}). 

In this picture, we also obtain an interpretation of the individual fibers of the 
meromorphic connection underlying the transform. Indeed, fix $\xi_0\in\Ct$, and let 
$(\xi-\xi_0)$ stand for the left ideal of ${\C}[\xi]\langle\partial_{\xi}\rangle$ 
generated by the monomial $\xi-\xi_0$. Define the fiber of $\Mt$ over $\xi_0$ by 
$$
   \Mt_{\xi_0}=\Mt/(\xi-\xi_0)\Mt. 
$$
Then, we have the isomorphism 
$$
   \Mt_{\xi_0}\cong \coker(M\xrightarrow{\partial_{x}-\xi_0}M). 
$$

\section{Meromorphic connections and $\D$-modules}\label{sec:merconn}
In this section we give the link between the objects studied in the previous sections. 
Let $(E,\nabla )$ be a holomorphic bundle endowed with a meromorphic integrable connection 
on $\CP$, with poles in $P\cup \{ \infty \}$ as described in Section \ref{sec:nahm}. 
Denote by $E(\ast P \cup \{ \infty \})$ the meromorphic bundle 
$E\otimes_{\Hol} \Hol(\ast P \cup \{ \infty \})$ with poles in $P\cup \{ \infty \}$. 
By definition, $\nabla$ acts on $E(\ast P \cup \{ \infty \})$. 
Denote by $\D_{\CP}$ the sheaf of holomorphic differential operators on $\CP$: it is 
the sheaf whose local sections are of the form $\sum_{i}a_{i}(t)\partial_{t}^{i}$, with 
$t$ a local holomorphic coordinate and $a_{i}$ holomorphic functions. 
$E(\ast P \cup \{ \infty \})$ then becomes a $\D_{\CP}$-module in the usual way: 
a holomorphic function $a(t)$ acting by multiplication, and $\partial_t$ acting by 
$\nabla_{\partial/\partial t}$. 
This $\D_{\CP}$-module is called the meromorphic bundle associated to $(E,\nabla )$, denoted $\DM$. 
It has the following properties: 
\begin{enumerate}
\item{It is a \emph{holonomic}  module, i.e. any local section is annihilated by 
   a local differential operator.}
\item{It has only regular singularities at finite distance, and an irregular singularity 
   of Poincar\'e rank $1$ at infinity.}
\end{enumerate}

Here is a fundamental notion that will underlie our construction: 
\begin{dfn}\label{extmindef}
The $\D_{\CP}$-module $\DN$ is said to be an \emph{extension} of $\DM$ if the meromorphic bundle 
$\Hol_{\CP}(\ast P \cup \{ \infty\}) \otimes_{\Hol_{\CP}}\DN$ induced by $\DN$ is equal to $\DM$. 
An extension $\DN$ is called \emph{minimal} if it admits no submodule or quotient module 
supported in a point. 
\end{dfn}
It is known (see e.g. Section 1 of \cite{sab-hi}) that a minimal extension always 
exists and is unique up to isomorphism; it will usually be denoted by $\DM_{\min}$. 

\subsection{Regular singularities}
We will now recall the local structure of meromorphic flat bundles %and their minimal extension 
and introduce the notion of parabolic structure for these objects, in the 
presence of a regular singularity. 
%For a reference on these topics, see \cite{sab}. 
In this section, we fix a $p\in P$, a small disk 
$\Delta \subset \C$ centered at $p$ and containing no other element of $P$, and 
we suppose that the meromorphic flat bundle $\DM$ has a regular singularity at $p$. 
$\Hol$ and $\D$ will stand for $\Hol_{\Delta}$ and $\D_{\Delta}$ respectively. 
\begin{dfn}
A \emph{lattice} of the meromorphic flat bundle $\DM|_{\Delta}$ is a coherent 
$\Hol$-module $F$ such that $F \otimes_{\Hol} \Hol(\ast \{ p\}) =\DM|_{\Delta}$. 
A lattice $F$ is called logarithmic if the connection $1$-form of $\nabla$ 
written in any local trivialization of $F$ has a first-order pole. 
\end{dfn}
If a logarithmic lattice $F$ is chosen for $\DM$, then there is a well-defined residue map of 
$\nabla$: 
$$
    \res(\nabla):F|_p \lra F|_p 
$$
induced by the action of $x \partial_x$. 
Modifying $F$ by a $\otimes \Hol(-m\{ p\})$ amounts to adding $m$ to all eigenvalues of $\res(\nabla)$. 
By a theorem of Deligne (see Corollary 2.21 \cite{sab2}), it is possible to choose a logarithmic lattice $F$ 
such that all eigenvalues of $\res(\nabla)$ lie in the interval $]-1,0]$. 
Hence, we obtain a finite set $\A_{j}$ of complex numbers with real part in $]-1,0]$, such that 
$\A_{j} + \Z$ is the set of eigenvalues of the residue of $\nabla$ on the various logarithmic 
lattices of $\DM$. For $a \in \A_{j}$ and $m\in \Z$, the generalized eigenspaces of 
$\res(\nabla)$ corresponding to the eigenvalues $a$ and $\alpha= a + m$ on the corresponding lattices
are isomorphic through the (possibly meromorphic) gauge transformation $(x-p)^{m}$. 
For any $\alpha \in \A_{j}+ \Z$, let $\psi_{p}^{\alpha}\DM$ stand for the generalized eigenspace
of $\res(\nabla)$ on the lattice $F\otimes\Hol(-m\{ p\})$ for the eigenvalue
$\alpha$, where $F$ is Deligne's lattice and $m=\lceil\Re\alpha\rceil$ 
is the smallest integer greater than or equal to the real part of $\alpha$. 
In other words, there exists an isomorphism 
$$
    (x-p)^{m} : \psi_{p}^{\alpha-m}\DM \lra \psi_{p}^{\alpha} \DM ,
$$
and $a=\alpha-m\in\A_{j}$. 
The assumption that the parabolic structure in $p$ is compatible with the connection means 
precisely that for all eigenvalues $\mu^{j}_{k}$ of $A^j$ it defines a decreasing exhaustive 
filtration $\psi_{p_{j}}^{\mu^{j}_{k},\bullet }\DM$ indexed by $\beta \in [0,1[$. 
By resonance-freeness, for all $\alpha\in\A_{j}+ \Z$ there exists a unique $\mu^{j}_{k}$ and 
$N\in\Z$ such that 
\begin{equation}\label{alphamuN}
   \alpha = \mu^{j}_{k} + N. 
\end{equation}
It follows that the parabolic filtrations on the $\mu^{j}_{k}$-eigenspaces extend to filtrations 
$\psi_{p}^{\alpha ,\bullet }$ indexed by real numbers in the interval $[N,N+1[$ for any 
$\alpha \in \A_{j}+ \Z$ by declaring that the above isomorphisms restrict to filtered isomorphisms 
\begin{equation}\label{filtrisom}
      (x-p)^{N} : \psi_{p}^{\mu^{j}_{k},\beta}\DM \lra \psi_{p}^{\alpha,\beta+N}\DM ,
\end{equation}
where $N$ is the integer appearing in (\ref{alphamuN}). 
\begin{dfn}
A \emph{parabolic structure} at a regular singular point $p$ of the meromorphic bundle $\DM$ is 
a collection for all $\alpha\in\C$ of a decreasing exhaustive filtration $\psi_{p_{j}}^{\alpha,\bullet }\DM$ 
of the generalised eigenspace $\psi_{p_{j}}^{\alpha}\DM$ of the residue indexed by elements $\beta$ 
of an interval of the form $[k,k+1[$ where $k$ is an integer, such that for all 
$\alpha\in\C,\beta\in\R$ and $N\in\Z$ the map 
$$
      (x-p)^{N} : \psi_{p}^{\alpha+N,\beta}\DM \lra \psi_{p}^{\alpha,\beta+N}\DM 
$$
be an isomorphism. 
\end{dfn}
Hence, we extended the parabolic structure of the holomorphic bundle $E$ as a parabolic structure of $\DM$. 
\begin{rem}
By admissibility of the parabolic structure of $E$ and because $\nabla$ is resonance-free, 
the filtration $\psi_{p}^{0,\bullet }$ is the trivial filtration 
(i.e. $\gr^0 \psi_{p}^{0,\bullet } = \psi_{p}^{0}$). 
\end{rem}

\subsection{Irregular singularity}

Here --- although it would be possible to work in a more general framework ---, 
we restrict to the  case of a Poincar\'e rank $1$ irregular singularity with 
diagonalizable polar part, as described in Section \ref{sec:nahm}. 

It is known that after taking the formal completion of the meromorphic bundle 
with respect to the $x^{-1}$-adic valuation 
(and possibly after a finite ramified cover, that we shall omit), 
the irregular connection can be put in the form (\ref{connirr}), 
without extra holomorphic terms.  
Let us choose such a trivialization of the completion of $\DM$. 
Then, in this basis the germ of connection  $\nabla- A \d z$ has regular singularity. 
Denote by $\widehat{\psi}^{\mu}_{\infty}$ the generalized eigenspace corresponding to 
the eigenvalue $\mu$ of this connection. The formal connection $\nabla- A \d z$ 
decomposes as direct sum of the $n'$ germs of connections with regular singularity 
$\d + C_{l}{\d z}/{z}$ for $l=1,\ldots , n'$ where $C_{l}$ is the diagonal matrix 
$\diag(\mu^{\infty}_{1+a_{l}},\ldots,\mu^{\infty}_{a_{1+l}})$ of size $a_{1+l}-a_{l}$. 
Since $\nabla$ is supposed to be resonance-free at infinity, none of the $\mu^{\infty}_k$ 
can be an integer. 
\begin{dfn}
A \emph{parabolic structure at infinity} on $\DM$ is the data of parabolic structures 
for the meromorphic bundles associated to all of the regular singular connections 
$\d + C_{l}{\d z}/{z}$ for $l=1,\ldots , n'$. 
%The parabolic structure at infinity is said to be admissible if all the 
%parabolic structures of the $\d + C_{l}{\d z}/{z}$ are admissible. 
\end{dfn}

\begin{rem}
Notice that the parabolic structures of the regular singular connections $\d + C_{l}{\d z}/{z}$ 
are then all admissibile, because by resonance-freeness at infinity the $0$-eigenvalues of their 
residues are zero, hence they clearly inject into the graded of the parabolic structure for the weight $0$. 
\end{rem}

\begin{lmm}
Near $\infty$ the minimal extension $\DM_{\min}$ agrees with the meromorphic bundle $\DM$. 
\end{lmm}
\begin{proof}
As none of the eigenvalues of the residue of $\nabla$ at $\infty$ are integers, 
the successive applications of the residue on any eigenvector yield sections 
with poles of arbitrarily high order. Hence, any $\D_{\CP}$-module containing 
$E$ necessarily contains $\DM$ as well. 
\end{proof}

An obvious consequence of the Lemma is that near $\infty$ the parabolic structure on the 
meromorphic bundle induces a parabolic structure on the minimal extension. 

\subsection{Stability}\label{ssec:stab}
Once a parabolic structure of a meromorphic integrable bundle $\DM$ in all singularities 
(including infinity) is given, it is possible to define its degree:
$$
    \deg(\DM ) = -\sum_{p \in P \cup \{\infty \}} \sum_{\beta \in [0,1[} 
    \sum_{\alpha \in \A_{p}+ \Z} \alpha \dim (\gr^{\beta} \psi_{p}^{\alpha}\DM)
$$
where $\gr^{\beta}\psi_{p}^{\alpha}\DM=\psi_{p}^{\alpha,\beta}\DM/\psi_{p}^{\alpha,>\beta}\DM$ is 
the graded vector space of the parabolic filtration $\psi_{p}^{\alpha,\bullet}$ 
corresponding to the weight $\beta$. Similarly, one can define the parabolic degree of $\DM$:
$$
    \pdeg(\DM ) = \deg(\DM ) + \sum_{p \in P \cup \{\infty \}} \sum_{\beta \in [0,1[} 
      \beta \dim (\gr^{\beta}_p \DM). 
$$
Finally, as usual, we define the (parabolic) slope of $\DM$ as the quotient of its parabolic degree 
and its rank. This then allows us to define the notion of stability: 
\begin{dfn}
A meromorphic integrable bundle $\DM$ endowed with a parabolic structure is said to be \emph{parabolically 
stable} if all nontrivial integrable subbundles $\DN$ have slope smaller than $\DM$. 
\end{dfn}

Here $\DN$ is endowed with the induced parabolic filtration: 
$$
   \psi_{p}^{\alpha,\beta}\DN = \psi_{p}^{\alpha}\DN \cap \psi_{p}^{\alpha,\beta}\DM. 
$$

\section{Minimal parabolic Laplace transform}\label{sec:parlapl}
Since $\DM$ is the meromorphic bundle associated to a meromorphic connection on a 
holomorphic bundle, and $\DM_{\min}$ its minimal extension, they are both modules 
over the sheaf $\D_{\CP}(\ast \infty)$ of meromorphic differential operators with 
a pole at infinity. The global sections of $\D_{\CP}(\ast \infty)$ form the Weyl algebra 
${\C}[x]\langle\partial_{x}\rangle$,
therefore the sets $M$ and $M_{\min}$ respectively of global sections on $\CP$ of $\DM$ and $\DM_{\min}$ 
admit a ${\C}[x]\langle\partial_{x}\rangle$-module structure. 

In Section \ref{sec:lapl} we defined Laplace transform $\Mt$ (respectively, inverse  Laplace transform 
$\check{N}$) for a ${\C}[x]\langle\partial_{x}\rangle$-module $M$ 
(respectively, a ${\C}[\xi]\langle\partial_{\xi}\rangle$-module $N$). 
\begin{dfn}
The \emph{minimal Laplace transform} of a meromorphic flat bundle $\DM$ is the meromorphic flat bundle 
$\widehat{\DM}$ associated to the ${\C}[\xi]\langle\partial_{\xi}\rangle$-module $\widehat{M_{\min}}$. 
Similarly, the inverse minimal Laplace transform of a meromorphic flat bundle $\DN$ 
on $\CPt$ is the meromorphic flat bundle $\check{\DN}$ associated to the 
${\C}[x]\langle\partial_{x}\rangle$-module $\check{N}_{\min}$. 
\end{dfn}

This is similar to N. Katz' Fourier transform in the $\ell$-adic case \cite{kat}. 

Let us now fix the parameters appearing in (\ref{aj}) and in the matrices $A,C$ of (\ref{connirr})
and introduce two categories: 
\begin{enumerate}
\item{$\mbox{HP}(\C)$, whose objects are holomorphic bundles with meromorphic connection 
      $(E,\nabla)$ over $\CP$, with resonance-free semi-simple logarithmic singularities in $P$
      given by (\ref{aj}), and a resonance-free irregular singularity of rank $1$ at infinity with 
      diagonalizable polar part given by (\ref{connirr}), endowed with an admissible parabolic structure 
      at all singular points}\label{def:hp}
\item{$\mbox{DMP}(\C)$, whose objects are meromorphic integrable bundles $\DM$ admitting 
      a lattice with the type of singularities as in $\mbox{HP}(\C)$, endowed with an admissible 
      parabolic structure in all singular points.}\label{def:dmp}
\end{enumerate}
Morphisms in $\mbox{HP}(\C)$ are holomorphic bundle maps compatible with connections and
parabolic structure, while 
in $\mbox{DMP}(\C)$, maps of $\D$-modules compatible with parabolic structure. 
Let us also introduce the full sub-categories 
$\mbox{HP}_{0}^{st}(\C)$, $\mbox{DMP}_{0}^{st}(\C)$ consisting of stable objects of parabolic degree $0$. 

\begin{prp}\label{prop:equiv}
The categories $\mbox{HP}(\C)$ and $\mbox{DMP}(\C)$ are equivalent.
The equivalence restricts to an equivalence between $\mbox{HP}_{0}^{st}(\C)$ and $\mbox{DMP}_{0}^{st}(\C)$. 
\end{prp}

\begin{proof}
Notice first that there exists an obvious functor 
\begin{equation}\label{functor}
    \mbox{HP}(\C)\to \mbox{DMP}(\C)
\end{equation}
mapping a holomorphic bundle endowed with a meromorphic connection and an admissible parabolic 
structure to the associated meromorphic flat bundle $\DM$ with parabolic structure as described 
in Section \ref{sec:merconn}. 
On the other hand, given an object $\DM$ of $\mbox{DMP}(\C)$, 
associate to it the integrable connection on a lattice as in its definition \ref{def:dmp}. 
To see that this is a quasi-inverse of (\ref{functor}) we need to show that such a lattice 
is uniquely determined. Near a logarithmic point $p_j$ this follows from the fact that 
because of the non-resonance condition, the logarithm of the local monodromy $\exp(2\pi iA^j)$ 
with given eigenvalues $(0,\ldots ,0,2\pi i\mu^{j}_{r_{j}+1},\ldots ,2\pi i\mu^{j}_{r})$ 
is uniquely determined. 
Near the irregular singularity at infinity we use a method originally due to Levelt 
(c.f. Section 2 of \cite{lev}). 
Namely, suppose we are given lattices $F_1$ and $F_2$ of $\DM$ in a neighborhood of $\infty$ 
with the property that in a basis of both of them the connection $D$ reads as in (\ref{connirr}), 
with the same matrices $A$ and $C$. 
Since both $F_1$ and $F_2$ are lattices of $\DM$, these two bases are then linked by a 
meromorphic gauge transformation $g$ that can be expanded into Laurent series 
$$
   g(t)=\sum_{b=N}^{\infty}g_bt^b
$$
for sufficiently small values of $t=x^{-1}$ and some $N\in\Z$. By assumption, we have 
$$
   g^{-1}(At^{-2}\d t + Ct^{-1}\d t+O(1)\d t) g - g^{-1}\d g = At^{-2}\d t + Ct^{-1}\d t+O(1)\d t, 
$$
or equivalently 
$$
   \d g = (At^{-2}\d t + Ct^{-1}\d t+O(1)\d t) g - g (At^{-2}\d t + Ct^{-1}\d t+O(1)\d t). 
$$
Inserting the Laurent series of $g$ into the latter equation yields the equations 
\begin{align*}
   Ag_N - g_NA & = 0 \\
   Cg_N - g_NC + Ag_{N+1} - g_{N+1}A & = Ng_N, 
\end{align*}
and so on. The first of these equations implies that $g_N$ is in the Levi subgroup 
$$
   L=\Gl(a_2-a_1,\C)\times \cdots \times \Gl(a_{1+n'}-a_{n'},\C)
$$ 
corresponding to the eigenspace-decomposition of $A$. 
In the second expression, the term $Ag_{N+1} - g_{N+1}A$ lies in the complementary of $L$, 
therefore the projection of $Cg_N - g_NC$ to $L$ must be equal to the right-hand side. 
Since the restriction of $C$ to any eigen-subspace of $A$ has by assumption distinct 
eigenvalues modulo the integers (except for multiple eigenvalues), it follows that the 
only integer eigenvalue of the adjoint action of $C$ restricted to $L$ is $0$. 
Therefore, the second equation implies $N=0$; in particular, $g$ is holomorphic. 
The same argument applied to $g^{-1}$ shows that $F_1=F_2$. 

It is clear that this equivalence preserves parabolic degree and stability. 
\end{proof}

\subsection{Parabolic transform}\label{ssec:partr}

The transform defined above is compatible with admissible parabolic structures: if an 
admissible parabolic structure is given on $\DM$, then we can define one on its 
(inverse) minimal Laplace transform $\DMt$ by the following construction. 

\begin{enumerate}
\item{In regular singularities: these are the eigenvalues $\{\xi_{1}, \ldots \xi_{n'} \}$ of
$A$. By formal stationary phase (Prop. V.3.6, \cite{sab2}) and the isomorphism 
${\psi}^{\alpha}\DM\cong{\psi}^{\alpha}\DM_{\min}$ for all $\alpha  \notin \Z$, 
there exists a formal isomorphism 
$$
     \widehat{\psi}^{\alpha}_{\infty}\DM \cong \bigoplus_{l=1}^{n'} {\psi}^{\alpha}_{\xi_{l}}\DMt, 
$$
where we recall that  $\widehat{\psi}^{\alpha}_{\infty}$ is the generalized eigenspace corresponding 
to the eigenvalue $\alpha$ of the ``regular part'' $\nabla - A \d z$ of the connection at infinity. 
The parabolic filtration on the left hand side therefore defines one on the right-hand side. 
Finally, on the $0$-eigenspaces of the residue in $\xi_{l}$, we put the trivial parabolic structure.}
\item{At infinity: similarly to the regular case, formal stationary phase applied to the inverse 
Laplace transform, for all $\alpha \notin \Z$ there exists a formal isomorphism 
$$
     \widehat{\psi}^{\alpha}_{\infty}\DMt \cong \bigoplus_{j=1}^{n} {\psi}^{\alpha}_{p_{j}}\DM, 
$$
and the parabolic filtrations on the direct summands on the right-hand side induce 
a parabolic filtration on the space on the left-hand side.} 
\end{enumerate}

The parabolic degree and stability condition corresponding to the structure 
of this definition behave well with respect to Laplace transform. The following 
statements are due to C. Sabbah, Proposition 4.1 and Theorem 5.1 of \cite{sab}; 
for sake of completeness we reproduce them here together with their proof. 
In what follows, we let $\DM$ denote a parabolically stable meromorphic integrable bundle 
of parabolic degree $0$. 

\begin{prp}\label{lem:minext}
The Laplace transform $\widehat{M_{\min}}$ is a minimal extension at all regular 
singular points $\xi_l$. In particular, parabolic minimal Laplace transform 
and inverse parabolic minimal Laplace transform are inverses of each other. 
\end{prp}
\begin{proof}
Suppose to the contrary that $\widehat{M_{\min}}$ is not a minimal extension at 
some $\xi_l$, so that it admits a submodule (or a quotient module) supported at $\xi_l$. 
By inverse Laplace transform, it follows that $M_{\min}$ admits a submodule (or a quotient) 
$N$ of the form $\d+\xi_l \d x$ on some trivial bundle. 
In particular, $N$ has a pole only at $\infty$, and the residue of the connection is $0$. 
However, by resonance-freeness at infinity $M$ cannot admit such a sub-module. 

Since $\widehat{M_{\min}}$ is a minimal extension, it is clearly the minimal 
extension of $\widehat{\DM}$, whose inverse parabolic minimal Laplace 
transform is therefore $M_{\min}$ endowed with its original parabolic structure. 
This implies the second claim. 
\end{proof}

\begin{rem}
The statement also holds without the resonance-freeness assumption, but the 
proof then makes use of stability: 
consider the meromorphic integrable bundle $\DN$ associated to $N$. 
Since the parabolic structure of $M_{\min}$ is admissible, the induced parabolic structure 
on $\mathcal{N}$ only has $0$ weight. Therefore, $\mathcal{N}$ is a subbundle of parabolic 
weight $0$ of the meromorphic bundle $\DM$ associated to $M_{\min}$, contradicting stability. 
\end{rem}

\begin{prp}
If $\DM$ is a parabolically stable meromorphic integrable bundle of parabolic 
degree $0$, then so is its parabolic minimal Laplace transform $\DMt$. 
\end{prp}

\begin{proof}
By the stationary phase formula for ordinary Laplace transform \cite{malgr}, the 
non-zero eigenvalues of the residue of the connection at the regular singularities 
are switched with the eigenvalues of the residue at infinity. 
In particular, the sum of the eigenvalues (counted with multiplicity) at all 
singularities is the same number for $\DM$ and $\DMt$. 
By our definition of parabolic structure for $\DMt$ and by admissibility, 
the same holds for the sum of all parabolic weights. Indeed, as the 
only weight on ${\psi}^{0}_{p_{j}}$ is $0$, the two total sums only differ 
in terms equal to $0$. 
In view of the formulae of Subsection \ref{ssec:stab}, 
this proves equality of the parabolic degree. 

It follows from Proposition \ref{lem:minext} that $\widehat{M_{\min}}$ is the minimal extension of $\DMt$. 
Let $\widehat{\DN}$ be a non-trivial meromorphic integrable subbundle of $\widehat{\DM}$. 
Then, its minimal extension $\widehat{\DN}_{\min}$ at all regular singular points is a 
non-trivial submodule of $\widehat{\DM}_{\min}$. Hence, the inverse Laplace transform 
$\DN_{\min}$ is a non-trivial submodule of $\DM_{\min}$. In particular, the corresponding 
meromorphic integrable bundle $\DN$ is a non-trivial subbundle of $\DM$. 
Clearly, $\DN$ is the inverse minimal Laplace transform of $\widehat{\DN}$. 

\begin{lmm}
The parabolic structure induced on $\DN$ from $\DM$ agrees with the 
inverse minimal Laplace transformed parabolic structure coming from $\widehat{\DN}$. 
\end{lmm}
\begin{proof}
We only give a proof of this fact under the extra assumption that for any $j\in\{1,\ldots ,n\}$ 
all the non-zero eigenvalues $\mu^j_k$ are distinct (not only modulo $\Z$); the general case 
is a consequence of functoriality of the stationary phase formula, c.f. Subsection \ref{ssec:prop4}. 

For any fixed $p_j\in P$, the set $\A_{j}(\DN)$ of eigenvalues of the residue 
of $\DN$ is contained in $\A_{j}(\DM)$, and for each $\alpha \in \A_{j}(\DN)+ \Z$, 
we have $\psi_{j}^{\alpha}\DN\subset\psi_{j}^{\alpha}\DM$. Since by assumption all of these 
eigenspaces are either $0$ or $1$-dimensional, it is clear that the only parabolic 
weight for the induced parabolic structure of a non-zero $\psi_{j}^{\alpha}\DN$ 
is equal to the parabolic weight of $\psi_{j}^{\alpha}\DM$, 
which in turn is equal to the parabolic weight of $\psi_{\ti}^{\alpha}\widehat{\DM}$ 
restricted to the $p_j$-eigenspace of the leading term of the residue. 
This latter, however, is the parabolic weight of 
$\psi_{\ti}^{\alpha}\widehat{\DN}$ restricted to the $p_j$-eigenspace of the 
leading term of the residue (notice that this latter space is non-zero if and only if 
$\psi_{j}^{\alpha}\DN\neq 0$ in view of the stationary phase formula).
\end{proof}

As above for $\DM$, the Lemma then implies equality between the parabolic 
degrees $\pdeg(\DN)$ and $\pdeg(\widehat{\DN})$. By stability of $\DM$, we have 
$\pdeg(\DN)<0$, hence also $\pdeg(\widehat{\DN})<0$, whence stability of $\widehat{\DM}$. 
\end{proof}

We are ready to state our main result. Consider the diagram of functors 

\begin{equation}\label{diagram}
  \xymatrix@C+3cm@R+2cm{ \mbox{HP}_{0}^{st}(\C) \ar@/_/[d] \ar@/^/[r]^{\mbox{Nahm}}
               & \mbox{HP}_{0}^{st}(\Ct) \ar@/^/[l]^{\mbox{inverse Nahm}}
               \ar@/_/[d] \\
              \mbox{DMP}_{0}^{st}(\C) \ar@/^/[r]^{\mbox{minimal Laplace}} \ar@/_/[u] 
              & \mbox{DMP}_{0}^{st}(\Ct) \ar@/^/[l]^{\mbox{inverse minimal Laplace}} \ar@/_/[u] 
  }
\end{equation} 
where the vertical arrows are the equivalences defined in Proposition \ref{prop:equiv}. 
\begin{thm}\label{thm:main}
Diagram (\ref{diagram}) commutes. 
\end{thm}
\begin{rem}
The coincidence at the level of ranks follows by comparing the rank formulae for 
Laplace (Proposition V.1.5. of \cite{malgr}) and 
Nahm transforms (formula (1.10) of \cite{Sz}).
\end{rem}

\section{Outline of the proof}
The main idea of the proof consists in comparing both constructions to the first 
hypercohomology of the twisted de Rham complex of $\DM_{\min}$. Set 
$$
   \mbox{DR}(\DM_{\min},\partial_x-\xi) = (\DM_{\min}\xrightarrow{\partial_x-\xi}\DM_{\min})
$$
for the de Rham complex of $\DM_{\min}$ twisted by $-\xi$, 
where the two non-zero sheaves are placed in degrees $0$ and $1$. 
For any sheaf complex $\mathcal{C}$ on $\CP$, write $\H^m(\mathcal{C})$ for its 
hypercohomology of degree $m$. On the other hand, for any $\xi \in \CPt$ set 
\begin{equation*}
    \nabla_{\xi} = \nabla-{\xi} \d x,
\end{equation*}
and call this operator the twisted integrable connection. 
Recall from Section \ref{sec:nahm} that we denote by $\Pt $ the set $\{ \xi_{1}, \ldots ,\xi_{n'} \}$ of 
regular singularities of the transformed object. 
Recall also that Nahm transform of a holomorphic bundle $(E,\nabla)$ 
with meromorphic connection on $\CP$ is a holomorphic bundle $(\Et,\Nt)$ with meromorphic 
connection on $\CPt$: the fiber over a fixed $\xi\in\Ct$ is the first $L^2$-cohomology 
for the Euclidean metric of the elliptic complex (\ref{ellcompl}), holomorphic structure is induced by the 
trivial holomorphic structure with respect to $\xi$, $\Nt$ is induced by $\hat{\d}^{1,0}-x\d\xi$, 
and the transformed metric is defined by the $L^2$-norm of the harmonic representative of an 
$L^2$-cohomology class. 
Recall finally that $(\Et,\Nt)$ admits a holomorphic extension induced by the harmonic metric $\hat{h}$: 
holomorphic sections of $\Et$ at a singular point are holomorphic sections in a neighborhood of 
that point whose norm is bounded. 
In particular, meromorphic sections of $\Et$ at a singular point are holomorphic sections in a 
neighborhood of that point whose norm  grows at most polynomially with respect to the inverse of 
the distance to the point. 
For $j=1,2$, denote by $\pi_{j}$ projection to the $j$-th factor in the product $\CP \times \CPt$. 

Our first aim is to define a natural holomorphic extension to $\CPt$ of the 
bundle of first hypercohomologies of $\mbox{DR}(\DM_{\min},\partial_x-\xi)$. 
Denote the sheaf of holomorphic sections of the vector bundle $E$ also by $E$, 
and let $F$ be the sheaf over $\CP$ whose local sections are meromorphic sections of 
$E \otimes \Om^{1}$ with at most a double pole at infinity and at most a simple pole 
in the points of $P$, such that the residue at a point $p_j \in P$ be in $\im (A^j)$. 
We then have a sheaf map 
\begin{equation}\label{nablaxi}
    \nabla_{\xi} : E \lra F.
\end{equation}
The bundle of first hypercohomologies of $\mbox{DR}(\DM_{\min},\partial_x-\xi)$ 
carries a natural holomorphic structure on $\Ct\setminus \Pt$ induced 
by the trivial holomorphic structure of $\DM_{\min}$ with respect to $\xi$. 
The same holds for the bundle of first hypercohomologies 
$\H^1(E \xrightarrow{\nabla_{\xi}} F)$. 

\begin{prp}\label{prop:1}
The complexes $\mbox{DR}(\DM_{\min},\partial_{x}-{\xi})$ and 
$E \xrightarrow{\nabla_{\xi}} F$ of sheaves on $\CP$ are quasi-isomorphic. 
Furthermore, the natural holomorphic structures of the bundles of first hypercohomologies 
over $\Ct\setminus\Pt$ agree under this isomorphism. 
\end{prp}

Via this proposition the holomorphic bundle of first hypercohomologies 
$$
\H^{1}(\mbox{DR}(\DM_{\min},\partial_{x}-{\xi}))
$$ 
on $\Ct \setminus \Pt$ admits a natural holomorphic extension to $\Ct$. 
In particular, this holomorphic extension will also induce a meromorphic 
extension over these points. 
At a regular singular point $\xi_{l}$ the holomorphic extension is defined as the 
sections whose value at $\xi_{l}$ is in the finite-dimensional vector space 
$\H^{1}(E \xrightarrow{\nabla_{\xi_{l}}} F)$. 
It remains to define a holomorphic extension at $\ti$. We follow Subsection
4.4.2 of \cite{Sz} and slightly modify the map (\ref{nablaxi}) near infinity 
as follows: let $U_0=\Ct$ and $U_{\ti}=\CPt\setminus \{0\}$ be the standard 
affine open charts of $\CPt$ with coordinates $\xi$ and $\zeta=\xi^{-1}$, 
and let $s_0$ and $s_{\ti}$ stand for the canonical global sections of the sheaf $\Hol_{\CPt}(\ti)$ 
satisfying 
\begin{align*}
    s_0(\xi)=\xi && s_{\ti}(\xi)=1 && \mbox{on } U_0\\
    s_0(\zeta)=1 && s_{\ti}(\zeta)=\zeta && \mbox{on } U_{\ti}. 
\end{align*}
On the pullback bundle $\pi_{1}^{\ast}E$ over $\CP \times \CPt$, introduce the sheaf map 
\begin{equation} \label{fibconoyau}
    \nabla_{\bullet} = \nabla \otimes s_{\ti} - \d z \otimes s_0: 
         \pi_{1}^{\ast}E \lra \pi_{1}^{\ast}F \otimes \Hol_{\CPt}(\ti). 
\end{equation}
Over the affine $U_0$ and  relative to $\CPt$ the map $\nabla_{\bullet}$ is a connection 
whose restriction to the fiber $\CP \times \{ \xi \}$ is $\nabla - \xi \d x$. 
\footnote{It might look disturbing that over the chart $U_{\ti}$ this is no longer 
a connection, but instead a $\zeta$-connection. 
Notice however that we are mainly interested in the first hypercohomology of this complex, 
and the hypercohomology of a map is independent of scalar rescaling.} 
We may then define the holomorphic extension of the bundle of hypercohomologies 
$\H^{1}(\mbox{DR}(\DM_{\min},\partial_{x}-{\xi}))$ to $\ti$ to consist of those holomorphic 
sections near $\ti$ that converge to a vector in 
$$
    \H^1(E\xrightarrow{\d z}F)\cong \oplus_j\im\res(\nabla,p_j), 
$$
the first hypercohomology of the fiber of (\ref{fibconoyau}) over $\ti$. 
This holomorphic extension is clearly induced by the holomorphic
extension of the original bundle $E$ at the points $p_j$. 

The next ingredient of the proof is an identification of the first $L^{2}$-cohomology 
$\Et_{\xi}=L^{2}H^{1}(D_{\xi})$ of (\ref{ellcompl}) for the Euclidean metric with the first 
hypercohomology of the twisted de Rham complex of $\DM_{\min}$. 

\begin{prp} \label{prop:2}
For $\xi \in \Ct \setminus \Pt$, the vector spaces $\Et_{\xi}=L^{2}H^{1}(D_{\xi})$ 
(for the Euclidean metric) and $\H^{1}(\mbox{DR}(\DM_{\min},\partial_x-\xi))$ are isomorphic. 
The corresponding holomorphic bundles over $\Ct \setminus \Pt$ are isomorphic. 
Finally, this isomorphism respects the natural meromorphic extensions to $\CPt$. 
\end{prp}

The last step in identifying the meromorphic bundles underlying the two transforms is: 

\begin{prp} \label{prop:3}
For all $\xi\notin\Pt$, 
the first hypercohomology $\H^{1}(\mbox{DR}(\DM_{\min},\partial_x-\xi))$ is isomorphic to the fiber
$(\widehat{M_{\min}})_{\xi}$ of the minimal Laplace transform in $\xi$. The corresponding 
holomorphic bundles over $\Ct \setminus \Pt$ are isomorphic. Finally, the natural meromorphic 
extensions over $\CPt$ of the underlying meromorphic bundles are also isomorphic. 
\end{prp}

Having thus obtained isomorphism of the meromorphic bundles, we move on to 
match the transformed connection with the $\partial_{\xi}$-action. 
Proposition 3.5 of \cite{Sz} states that the transformed connection $\Nt$ is induced 
on $L^{2}$-cohomology by taking the quotient of $\dt - x \d \xi$. The isomorphism 
of Proposition \ref{prop:2} maps it to the connection induced by taking the quotient of 
$\dt - x \d \xi$ in hypercohomology. 
This, coupled with the interpretation of Laplace transform as a cokernel of
Subsection \ref{conoyauinterpr}, shows that  
over the open set $\Ct \setminus \Pt$ the connection of the Nahm transform is mapped 
to the $\partial_{\xi}$-action of the minimal Laplace transform. 

To see that the square of Theorem \ref{thm:main} commutes, 
the last thing to show is that the minimal 
extension over the singularities of the meromorphic bundle with connection
associated to $(\Et,\Nt )$ is isomorphic to $\widehat{M_{\min}}$. 
Since the underlying meromorphic connections agree and by uniqueness up to isomorphism of 
the minimal extension of a meromorphic bundle with connection, Proposition \ref{lem:minext} 
implies the claim. 

Finally, it remains to compare the parabolic structures of the transformed objects: 

\begin{prp} \label{prop:4}
The parabolic structure of $(\Et ,\Nt )$ is mapped under this correspondence into the one
defined in Subsection \ref{ssec:partr}. 
\end{prp}

\section{Proof of the propositions}

We will carry out all proofs for the direct transforms; the inverse transforms 
are known to agree with the direct ones up to a sign $\xi\leftrightarrow -\xi$, 
therefore the same arguments imply the statements for the inverse transforms. 

\subsection{Proposition \ref{prop:1}}

First, let us define a sequence of sub-sheaves of $\DM_{\min}$: for $m \geq 0$ set $F_{m}$ for 
the sheaf of sections of $\DM_{\min}$ with pole of order at most $m$ at all singular points. 
These sheaves define an exhaustive filtration of $\DM_{\min}$ : 
$$
     E=F_{0}\subset F_{1}\subset F_{2}\subset \ldots \subset \DM_{\min} 
$$
such that $\partial_{x} -\xi$ maps $F_{m}$ into $F_{m+1}$ for all $\xi\notin\Pt$. 
Notice that by definition $F=F_{1}\otimes\Om^{1}$. 
Moreover, the quotient sheaves $F_{m+1}/F_{m}$ are 
supported in the set of  regular singular points $P$, and the stalk at $p_j$ is of dimension $r-r_j$. 
Denote by $[\partial_{x} -\xi]$ the map induced by $\partial_{x} -\xi$ on these quotients. 
\begin{lmm}\label{ujlem:1.1}
For all $m>0$ the map 
$$
   [\partial_{x} -\xi]: F_{m}/F_{m-1} \lra F_{m+1}/F_{m}
$$
is an isomorphism of complex vector spaces. 
\end{lmm}
\begin{proof}
Since twisting by $\xi$ does not modify the principal term of $\partial_{x}$, it is sufficient 
to show the result for $\xi=0$. 
The stalk of $F_{m}/F_{m-1}$ at $p \in P$ is isomorphic to the image of the residue of $\nabla$ in 
$p$ (up to the $1$-form $\d z$), and the same thing holds for $F_{m+1}/F_{m}$. 
The action of $[\partial_{x}]$ after this identification is the matrix $\res(\nabla,p)-m\cdot\mbox{Id}$. 
Since $\nabla$ is resonance-free, this matrix admits no $0$ eigenvalue. Since
the source and target spaces are of the same dimension, it is an isomorphism. 
\end{proof}

Clearly, the complexes $E\xrightarrow{\nabla_{\xi}}F$ and $E\xrightarrow{\partial_{x}-{\xi}}F_{1}$ 
are quasi-isomorphic. 
By the lemma, for all $m>0$ the complexes $F_{m}\xrightarrow{\partial_{x}-{\xi}}F_{m+1}$ and 
$F_{m-1}\xrightarrow{\partial_{x}-{\xi}}F_{m}$ are quasi-isomorphic. Since $\DM_{\min} = \cup_{m} F_{m}$, 
the complex $\DM_{\min} \xrightarrow{\partial_{x}-{\xi}} \DM_{\min}$ is the inductive limit of its 
sub-complexes $F_{m}\xrightarrow{\partial_{x}-{\xi}}F_{m+1}$. The sequence of these latter being 
stationary (up to quasi-isomorphism) for $m\geq 0$, we conclude that  
$\DM_{\min} \xrightarrow{\partial_{x}-{\xi}} \DM_{\min}$ is quasi-isomorphic to 
$F_0\xrightarrow{\partial_{x}-{\xi}}F_{1}$, hence to $E\xrightarrow{\nabla_{\xi}}F$. 
This shows the fiber-wise statement. The statement about holomorphic structures is clear. 

\subsection{Proposition \ref{prop:2}}\label{ssec:prop:2}

Consider the double complex 
\begin{equation} \xymatrix{ {L}^{2}(\Omega^{0,1}\otimes E) \ar[r]^{\nabla_{\xi } }    &  {L}^{2}(\Omega^{1,1}\otimes E) \\
              {L}^{2}(E) \ar[r]^(.4){\nabla_{\xi }} \ar[u]^{\bar{\partial}^{E} } &
              {L}^{2}(\Omega^{1,0}\otimes E)  \ar[u]_{\bar{\partial}^{E} } \\
              {E} \ar[r]^{\nabla_{\xi }} \ar@{^{(}->}[u] & {F} \ar@{^{(}->}[u] .
} \label{compldouble}
\end{equation}
Here we use the convention that whenever ${L}^{2}$ is written at a particular place of a 
diagram, it is meant to be the intersection of the domains of the operators corresponding 
to all arrows that point out from that place. 

To prove equality between $L^2$-cohomology and hypercohomology it is sufficient to show two statements: 
\begin{lmm} \label{lem2.1}
For all $\xi \notin \Pt$ the vertical arrows in (\ref{compldouble}) are sheaf resolutions. 
\end{lmm}

\begin{rem}
A similar result is shown in Theorem 6.2 of \cite{z} for the Poincar\'e-metric near the punctures and 
with respect to a Hermitian metric with zero parabolic weight but non-trivial weight-structure, 
see also Theorem 1.2 of \cite{sab-hi}. However, we included a proof because we see no direct 
way of applying those results here. 
\end{rem}

\begin{lmm} \label{lem2.2}
The first cohomology of the simple complex associated to (\ref{compldouble}) is equal to 
$L^{2}H^{1}(D_{\xi})$. 
\end{lmm}

\begin{proof}[Proof of Lemma \ref{lem2.1}]
This is analogous to Lemma 4.12 of \cite{Sz}, where the same statement is proved for Higgs bundles. 
As the problem is local, and away from singularities all spaces involved are usual $L^2$-spaces 
(so exactness is guaranteed by usual $L^2$-theory), we focus on a neighborhood of a singular point
that we may suppose to be a disk. 

Let us first treat the column on the left, near a regular singular point $p=p_{j}$. 
Fix a holomorphic local trivialization $e_{1},\ldots e_{r}$ of $E$ 
in which the connection can be written $\nabla =\d + A^{j} {\d x}/{(x-p_{j})}$ up to 
holomorphic terms, where $A_{j}$ is the matrix in (\ref{aj}). In all that follows, we will 
omit the index $j$ of the singularity. Let $e = \sum_{k}\varphi_{k}e_{k}$ be a local 
$L^{2}$ section of $E$ in the kernel of $\dbar^{E}$: this means that the $\varphi_{k}$ 
are meromorphic functions. Recall that the $L^{2}$ condition  is measured with respect 
to the metric $h$, and this latter has asymptotic behavior bounded with 
$diag(|x-p|^{2\beta_{k}})$, where by admissibility $0 \leq \beta_{k} <1$ for $r_{j}<k$ 
and $ \beta_{k}=0$ otherwise. Hence the condition $e \in L^{2}$ is equivalent to the conditions 
\begin{itemize}
\item for $k\leq r_j$ the coefficients $\varphi_{k}$ are holomorphic at $p$; 
\item for $r_{j}<k$ such that $\beta_k=0$ the coefficients $\varphi_{k}$ are holomorphic at $p$; and 
\item for $r_{j}<k$ such that $\beta_k>0$ the coefficients $\varphi_{k}$ have at most a simple pole at $p$. 
\end{itemize}
Let us show that the requirement $\nabla_{\xi}e \in L^{2}$ implies even in the last case that $\varphi_{k}$ 
must be holomorphic; for this purpose consider their Laurent-series: 
$$
     \varphi_{k}(z) = \varphi_{k,-1}\frac{1}{x-p_{j}} + \varphi_{k,0} + \varphi_{k,1} (x-p_{j}) + \ldots
$$
By the local form of $\nabla$ in the trivialization $e_{1},\ldots e_{r}$, 
the coefficient of $e_{k}\d z$ for $r_{j}<k$ in $\nabla_{\xi}e$ is
\begin{align*}
     \left( \frac{\partial}{\partial x}+\frac{\mu_{k}}{x-p} \right)\varphi_{k}(x)+\sum_{l=1}^{r}h_{l}(x)\varphi_{l}(x), 
\end{align*}
where the $h_{l}$ are holomorphic functions. Plugging the Laurent-expansion into this yields 
$$
    (\mu_{k}-1)\frac{\varphi_{k,-1}}{(x-p)^{2}} +  O\left( \frac{1}{x-p} \right), 
$$
where $O(\cdot)$ denotes terms with at most a simple pole, depending on $\varphi_{l,s}$, 
$\mu_{l}$ and $h_{l}$. Such a $1$-form is in $L^2$ if and only if $(\mu_{k}-1)\varphi_{k,-1}=0$. 
Since $\nabla$ is resonance-free, this is equivalent to $\varphi_{k,-1}=0$, in different terms 
holomorphicity of $\varphi_{k}$. On the other hand, for $k\leq r_j$, taking into account 
that all $\varphi_{l}$ are holomorphic and that $\mu_k=0$, the coefficient of $e_{k}\d x$ 
in $\nabla_{\xi}e$ is clearly holomorphic. 
%$$
%    \varphi_{k,1} + O(x-p),
%$$
%in particular it is holomorphic. The sum of finitely many such expressions is again holomorphic. 
This proves that $e$ is in the $L^2$-domain of $\nabla$ if and only if it is holomorphic, 
hence the left column is exact in the term $L^{2}(E)$. 

Let us now study the second term of this column: given a section $f=\sum_{k}\phi_{k}e_{k}\d \bar{x}$ 
of $L^{2}(\Omega^{0,1}\otimes E)$, we would like to show that there exists an $L^2$ section 
$e = \sum_{k}\varphi_{k}e_{k}$ of $E$ such that $\nabla_{\xi}e \in L^{2}$ and $\dbar^{E}e=f$. 
As $e_{1},\ldots,e_{r}$ is a holomorphic trivialization (in particular, $\dbar^{E}$ is 
diagonal), we may suppose that $f= \phi e_{k}\d \bar{x}$ for some $k$. 
For $k\leq r_{j}$, by admissibility both $\nabla$ and the Hermitian metric are equivalent to 
the usual metric, so the usual $L^{2}$ Dolbeault lemma gives the result. 
For $r_{j}<k$, the section $\phi e_{k}$ is in $L^{2}$ with respect to $h$ if and only if 
$\phi |x-p|^{\beta_{k}}$ is an $L^{2}$ function. 
In the case $\beta_{k}>0$ the parabolic $L^2$ Dolbeault lemma (Claim 4.13 of \cite{Sz}) 
implies that there exists $\varphi$ such that $\dbar \varphi=\phi$ and 
$\varphi|x-p|^{\beta_{k}-1}\in L^{2}$. 
Finally, in the case $r_{j}<k$ and $\beta_k=0$ the result again follows from the usual 
$L^{2}$ Dolbeault lemma: for $f= \phi e_{k}\d \bar{x}\in L^2$ we find $\varphi\in L^2$ 
such that $\dbar \varphi=\phi$; now as by assumption $\nabla(f)$ is also in $L^2$ 
and the diagram (\ref{compldouble}) commutes we have that $\nabla\varphi$ is (up to subtracting 
a holomorphic $1$-form) also in $L^2$; whence surjectivity of $\dbar^{E}$. 

This finishes the proof of exactness of the left column. 

All the other statements contained in the lemma (i.e. that the right column is also 
exact in the regular singularities and that both columns are exact near infinity) are 
immediate modifications of the same arguments. We leave details to the reader. 
\end{proof}

\begin{proof}[Proof of Lemma \ref{lem2.2}]
The corresponding statement for Higgs bundles is Proposition 4.7 of \cite{Sz}.

As the differential of the simple complex associated to (\ref{compldouble}) is $D_{\xi}$, 
this is purely a question of domains. Denote by $\Comp_{\xi}$ the simple complex  associated 
to (\ref{compldouble}). An element of the degree $1$ term $\Comp_{\xi}^{1}$ of $\Comp_{\xi}$ 
is represented by a sum  
$$
f^{0,1}\d \bar{x} + f^{1,0}\d x \in L^{2}(\Omega^{0,1}\otimes E)\oplus L^{2}(\Omega^{1,0}\otimes E)
$$ 
such that 
$\nabla_{\xi}f^{0,1}\d \bar{x}\in L^{2}(\Omega^{1,1}\otimes E)$ and 
$\dbar^{E}f^{1,0}\d x\in L^{2}(\Omega^{1,1}\otimes E)$. 
Obviously, one then also has 
$D_{\xi} (f^{0,1}\d \bar{x}+ f^{1,0}\d x)\in L^{2}(\Omega^{1,1}\otimes E)$; in other words, 
there exists a tautological injection 
$$
    \Comp_{\xi}^{1} \lra L^{2}(D_{\xi})^{1}, 
$$
the latter space being the $L^{2}$ domain of $D_{\xi}$ on $1$-forms. 
Since terms of degree $0$ of the 
complexes $\Comp_{\xi}$ and $L^{2}(D_{\xi})$ are the same together with differentials, 
this injection induces an injection of degree $1$ cohomology spaces 
$$
    H^{1}(\Comp_{\xi}) \lra L^{2}H^{1}(D_{\xi}). 
$$
Let us show that this map is surjective too: for this, suppose $g=g^{0,1}\d \bar{x} + g^{1,0}\d x$ 
represents a cohomology class in $L^{2}H^{1}(D_{\xi})$; by definition, this means that $g$
is in $L^{2}$ and in the kernel of $D_{\xi}$. We would like to find a section 
$f=f^{0,1}\d \bar{x} + f^{1,0}\d x$ representing the same class, such that we have in addition 
$\nabla_{\xi}f^{0,1}\d \bar{x}\in L^{2}(\Omega^{1,1}\otimes E)$ and 
$\dbar^{E}f^{1,0}\d x\in L^{2}(\Omega^{1,1}\otimes E)$. Since this is guaranteed away from 
singularities by usual $L^2$-theory and because $\nabla_{\xi}$ is an isomorphism near infinity, 
it is sufficient to focus on a neighborhood of a regular singularity. 
In such a neighborhood, Lemma \ref{lem2.1} allows us to find $\varphi\in L^{2}(E)$ 
such that $\dbar^{E}\varphi = g^{0,1}\d \bar{x}$ and $\nabla_{\xi}\varphi\in L^{2}(\Omega^{1,0}\otimes E)$. 
The local $L^2$ section 
$$
f=(g^{0,1}\d \bar{x} - \dbar^{E}\varphi)+ (g^{1,0}\d x -\nabla_{\xi}\varphi)
$$
obviously satisfies the required conditions. 
\end{proof}

Let us come to the second statement of the proposition. 
Proposition 3.5 of \cite{Sz} gives the holomorphic structure of $\Et$: on the 
open set $\Ct \setminus \Pt$, it is the quotient in $L^{2}$-cohomology of the trivial 
holomorphic structure relative to $\CPt$ of $\pi_1^*E$. On the other hand, we defined 
the holomorphic structure on the bundle of first hypercohomologies 
$\H^{1}(E \xrightarrow{\nabla_{\xi}}F)$ as the quotient of the trivial 
holomorphic structure on $F$ relative to $\CPt$. This yields equality of the holomorphic 
structures away from singularities. 

Finally, in order to prove that the meromorphic extensions agree, it is sufficient to show 
that the $L^{2}$-norms of harmonic $1$-forms whose limit as $\xi \ra \xi_{l}$ (respectively 
$\xi \ra \ti$) is a vector of the space $\H^{1}(E\xrightarrow{\nabla_{\xi_{l}}}F)$, 
is of moderate growth. Indeed, in Section \ref{sec:nahm} the holomorphic extension of the 
Nahm transformed bundle $\Et$ to the singularities is defined by the condition that 
the $L^2$-norms of harmonic representatives be bounded; hence, meromorphic sections 
near a singular point are exactly the ones having at most polynomial growth. 
We only treat the case of the singularity at infinity, the case of the regular 
singularities being similar. 

The crucial point is that since the Laplace operator (\ref{Laplacian}) is (up to a constant $2$) 
well-known  to be equal to 
\begin{equation}\label{Laplacian''}
    \Delta_{\xi}'' = D_{\xi}''(D_{\xi}'')^* + (D_{\xi}'')^* D_{\xi}'' 
\end{equation}
where $D_{\xi}''$ is the operator defined in (\ref{D''}), 
the space of $L^2$ harmonic $1$-forms for the two operators are the same. 
In particular, for fixed $\xi$ the $1$-form $\hat{\sigma}^{\infty}_{k}(\xi,x)$ 
constructed in Section \ref{sec:nahm} is also a vector of $\Et_{\xi}$. However, as 
$\xi$ varies, these vectors are holomorphic for the Dolbeault holomorphic 
structure, and not the de Rham holomorphic structure, see Section 4.4 of \cite{Sz}. 
Recall that $\sigma^{j}_{k}$ was chosen to be a holomorphic section of 
$\mathcal{E}(p_j)$ near $p_j$ such that for each $\xi\in\Ct\setminus \Pt$ 
the element $\sigma^{j}_{k}(q_k(\xi))$ represent a non-vanishing class in 
$\coker(\theta(q_k(\xi)) - \xi \d x/2)$, and $q_k(\xi)$ depends analytically on $\xi$. 
Let us now modify the construction of $\hat{\sigma}^{\infty}_{k}(\xi,x)$ 
to obtain a smooth local section $\hat{\tau}^{\infty}_{k}(\xi,x)$ of $\Et$ in 
a punctured disk near $\ti$ in the following way: instead of $\sigma^{j}_{k}$, 
pick a local section $\tau^{j}_{k}$, holomorphic with respect to the de 
Rham holomorphic structure $\dbar^E$ near $p_j$, extending an eigenvector of 
the image of the residue of $\nabla$ at $p_j$. Notice that by resonance-freeness 
the images of $\res(\nabla,p_j)$ and of $\res(\theta,p_j)$ are isomorphic, 
and formula (1.8) of \cite{biqboa} means that a possible choice for $\tau^{j}_{k}$ is 
\begin{equation}\label{tausigma}
   \tau^{j}_{k}(x)\approx |x-p_j|^{\beta^j_k-\alpha^j_k}\sigma^{j}_{k}(x)
\end{equation}
where $\alpha^j_k$ is the fractional part of $\Re\mu^j_k$. In particular, 
for each $\xi$ we still have that $\tau^{j}_{k}(q_k(\xi))$ represents a 
non-trivial class in $\coker(\theta(q_k(\xi)) - \xi \d x/2)$. 
Therefore, the definitions of $v^j_k(\xi)$, $t^j_k(\xi)$ and $\hat{\sigma}^{\infty}_{k}(\xi,x)$ in 
Section \ref{sec:nahm} carry over to this situation to define $w^j_k(\xi)$, $s^j_k(\xi)$ 
satisfying the same equation 
$$
    D_{\xi}''(w^j_k(\xi,x)\d x + s^j_k(\xi,x)\d \bar{x})=0, 
$$
and therefore a local smooth section $\hat{\tau}^{\infty}_{k}(\xi,x)$ of $\Et$ near $\ti$. 
\begin{lmm}\label{lem:ext}
On $\Ct\setminus \Pt$, the local smooth section $\hat{\tau}^{\infty}_{k}(\xi,x)$ 
of $\Et$ is holomorphic in $\xi$ with respect to the de Rham holomorphic
structure. In particular, the set of local sections $\hat{\tau}^{\infty}_{k}(\xi,x)$ 
for all $p_j\in P$ and a basis $\tau^{j}_{k}$  of $\im(\res(\nabla,p_j))$ 
defines a holomorphic extension $\widehat{F}$ of $\Et|_{\Ct\setminus \Pt}$ at $\ti$. 
\end{lmm}
\begin{proof}
In the construction of $\hat{\tau}^{\infty}_{k}(\xi,x)$, the dependence 
on $\xi$ lies in the spectral point $q_k(\xi)$. 
\footnote{A subtle point here is that in (\ref{trsection}) the argument of the Gaussian function $\chi$ 
contains $|\xi|$, which is not holomorphic in $\xi$. However, locally near any $\xi_0\notin \Pt$ 
this factor could be kept constant $|\xi_0|$ independently of $\xi$, and the sections 
so obtained would be cohomologous to the sections (\ref{trsection}), so the non-holomorphic 
dependence is actually only apparent. The point in the choice of the argument of $\chi$ 
is just that its support be sufficiently separated from the other spectral points. 
Of course, the sections with $|\xi|$ replaced by $|\xi_0|$ do not extend to a neighborhood 
of the singular points $\Pt$, so for the construction of asymptotically harmonic sections 
only the choice in (\ref{trsection}) is suitable.} 
As this latter varies 
holomorphically with $\xi$, the monad construction of \cite{don-kro}, 
Subsection 3.1.3, shows that $\hat{\tau}^{\infty}_{k}(\xi,x)$ is 
holomorphic with respect to the holomorphic structure on $\Lisse(\hat{E})$ 
induced by the pull-back to $(\C\setminus P)\times(\Ct\setminus \widehat{P})$ 
of $E$, which is by definition the de Rham holomorphic structure of $\Et$. 

By Theorem 1.6 of \cite{Sz}, the rank of the transformed integrable bundle is equal to 
$\sum_{j=1}^n\dim(\im(\res(\nabla,p_j)))$, which is the number of the local sections 
$\hat{\tau}^{\infty}_{k}(\xi,x)$ for all $p_j\in P$ and $\tau^{j}_{k}\in\im(\res(\nabla,p_j))$. 
To see that these latter define a holomorphic extension, it is therefore sufficient to prove that they 
are linearly independent over $\C(\xi)$. This will follow from the fact that as $\xi\to\ti$ they are 
linked to the linearly independent sections $\hat{\sigma}^{\infty}_{k}(\xi,x)$ via (\ref{taukalapsigma}). 
\end{proof}

\begin{lmm}
The sections $\hat{\tau}^{\infty}_{k}(\xi,x)$ are asymptotically harmonic. 
\end{lmm}
\begin{proof}
Identical to the proof of Lemma 4.34 of \cite{Sz}, for this latter does not rely at all on 
the choice of vectors representing classes in $\coker(\theta(q_k(\xi)) - \xi \d x/2)$ as in 
(\ref{tausigma}). 
\end{proof}
This now implies that the growth of the $L^2$-norm of these 
sections is at most polynomial. Namely, we have: 
\begin{lmm}\label{lem:parwght}
The $L^2$-norm of $\hat{\tau}^{\infty}_{k}(\xi,x)$ over $\C$ with respect to 
the harmonic metric $h$ satisfies 
$$
    c|\xi|^{2-2\beta^j_k}\leq \Vert
    \hat{\tau}^{\infty}_{k}(\xi,x)\Vert^2_{L^2}\leq 
    C|\xi|^{2-2\beta^j_k}
$$
for some constants $0<c<C$. 
\end{lmm}
\begin{proof}
Similar to the proof of Theorem 4.32 of \cite{Sz}. The difference with that case 
lies only in the behavior of the metric $h$ on the sections $\tau^j_k$ and $\sigma^j_k$. 
In view of (\ref{tausigma}) and the asymptotic $q_{k}(\xi)-p_j\approx 2\lambda^{j}_{k}\xi^{-1}$ 
(Claim 4.27, loc. cit.) and because the construction of $w^j_k,s^j_k$ depends $\C$-linearly 
with the chosen cokernel vector $\tau^j_k$, the sections $\hat{\tau}^{\infty}_{k}(\xi)$ satisfy 
\begin{equation}\label{taukalapsigma}
   \hat{\tau}^{\infty}_{k}(\xi) \approx c_0 |\xi|^{-\beta^j_k+\alpha^j_k}\hat{\sigma}^{\infty}_{k}(\xi)
\end{equation}
as $\xi\to\ti$ for some constant $c_0\neq 0$ depending on $\lambda^{j}_{k},\alpha^j_k,\beta^j_k$. 
By Theorem 4.32 of loc. cit., we have 
$$
    c_1|\xi|^{2-2\alpha^j_k}\leq \Vert
    \hat{\sigma}^{\infty}_{k}(\xi,x)\Vert^2_{L^2}\leq 
    C_1|\xi|^{2-2\alpha^j_k}. 
$$
Comparing these estimates yields the desired asymptotic. 
\end{proof}

This finishes the proof of the proposition.

\subsection{Proposition \ref{prop:3}}

The hypercohomology long exact sequence gives a decomposition 
\begin{align*}
     \H^{1}(\DM_{\min} \xrightarrow{\partial_{x}-{\xi}} \DM_{\min}) \isom &
     \coker(H^{0}(\DM_{\min})\xrightarrow{\partial_{x} -\xi}H^{0}(\DM_{\min})) \oplus \\
     & \ker(H^{1}(\DM_{\min})\xrightarrow{\partial_{x} -\xi}H^{1}(\DM_{\min})). 
\end{align*}
By definition, $H^{0}(\DM_{\min})$ is $M_{\min}$, and resonance-freeness at infinity implies that 
in a neighborhood of infinity the equality $\DM_{\min}=\DM_{\min}(\ast \{ \infty \})$ holds, 
i.e. locally at infinity $\DM_{\min}$ coincides with the meromorphic bundle. In particular,
the first cohomology $H^{1}(\DM_{\min})$ vanishes, and it follows that 
$$
    \H^{1}(\DM_{\min} \xrightarrow{\partial_{x}-{\xi}} \DM_{\min}) \isom 
     \coker(M_{\min}\xrightarrow{\partial_{x} -\xi}M_{\min}).
$$
This gives the fiber-wise statement of the proposition. 

The holomorphic structure induced on the bundle $\H^{1}(E \xrightarrow{\nabla_{\bullet}}F)$ by 
taking the quotient of the trivial holomorphic structure of $\pi_{1}^{\ast}F$ relative to $\CP$ 
is mapped by the above isomorphism to the one induced by the trivial $\C[\xi]$-module structure 
of $\M_{\min}=M_{\min}\otimes_{\C}\C[\xi]$ on the cokernel of $\partial_{x} -\xi$. 
This latter is by definition the holomorphic structure of $\widehat{M_{\min}}$ on $\Ct \setminus \Pt$. 

It remains to match the meromorphic extensions over the singularities. 
For this purpose, we need to study what happens to the isomorphism of
Proposition \ref{prop:1} as $\xi \ra \xi_l$ for some $\xi_l \in \widehat{P}$
or as $\xi\to\ti$. Near the points $\xi_l$, it is known that $\DMt$ admits 
a regular singularity. On the other hand, by Theorem 1.6 of \cite{Sz} 
the connection $\widehat{\nabla}$ on the transformed extension of $\Et$ over $\xi_l$ 
is logarithmic; in particular, the meromorphic flat bundle $\Et(*\Pt)$ has 
regular singularities at all points of $\Pt$. By uniqueness of the meromorphic 
extension with regular singularity, the meromorphic flat bundles are 
naturally isomorphic over $\Ct$. 

Therefore, we only need to treat the case of the irregular singular point $\ti$.
For this purpose, consider the holomorphic bundle $E$ over $\CP$ as a 
$\C[x]$-submodule of $M_{\min}$. The lattice $\widehat{F}$ near $\ti$ introduced in 
Lemma \ref{lem:ext} can be described in the cokernel-interpretation of Subsection 
\ref{conoyauinterpr} as the lattice generated by the classes of $F\subset M_{\min}$, 
i.e. of sections admitting a pole of order $1$ at the points of $P$ 
with respect to the lattice $E$. 
By resonance-freeness, $M_{\min}$ is spanned as a ${\C}[x]\langle\partial_{x}\rangle$-module 
by $F$. Therefore, by Theorem V.2.7 of \cite{sab2} $\widehat{F}$ is a lattice 
of the Laplace transform $\widehat{M_{\min}}$ of $M_{\min}$ near infinity; 
this concludes equality of the meromorphic  extensions.

\subsection{Proposition \ref{prop:4}}\label{ssec:prop4}

This follows from comparing the construction of Subsection \ref{ssec:partr} with 
Subsection \ref{ssec:prop:2}. Notice first that it is sufficient to treat the 
case of the singularity at $\ti$. Indeed, the result then also follows for the 
finitely located singularities using inverse transforms. 

First, let us prove the stationary phase formulae used in the construction of 
Subsection \ref{ssec:partr}. As we have already obtained the isomorphism of 
Nahm and minimal Laplace transforms (without parabolic structure), it is 
sufficient to show the results for Nahm transform. These latter, as we will see, 
are a simple consequence of the stationary phase formulae for Nahm transform 
of parabolic Higgs bundles, in view of the non-Abelian Hodge theory of \cite{biqboa}. 
Hence, we need first to describe how one obtains the stationary phase formulae for 
Nahm transform of parabolic Higgs bundles. 
For this purpose, notice that by Claim 4.27 of \cite{Sz} all of the spectral points $q_m(\xi)$ 
converge to one of the $p_j$'s  as $\xi\to\ti$, and for fixed $1\leq j \leq n$ there are exactly 
$\rank(A^j)=r-r_j$ indices $m$ for which the corresponding $q_m(\xi)$ converges to $p_j$. 
It follows that the sections $\hat{\sigma}^{\infty}_{k}(\xi)$ introduced in Section 
\ref{sec:nahm}, for all $p_j\in P$ and a diagonalizing basis $\sigma^{j}_{k}$ of 
$\im(\res(\theta,p_j))$ define a lattice $\widehat{\mathcal{F}}$  at $\ti$ of the transformed 
Higgs bundle. Furthermore, by Theorem 4.31 of \cite{Sz} with respect to this lattice 
the polar part at infinity of the transformed Higgs field reads 
\begin{equation}\label{stphDol}
   -\frac{1}{2} \begin{pmatrix}
                p_{1} & & & & & & \\
                & \ddots & & & & & \\
                & & p_{1}  & & & & \\
                & & & \ddots & & & \\
                & & & & p_{n} & & \\
                & & & & & \ddots & \\
                & & & & & & p_{n}
                        \end{pmatrix} \d \xi - 
               \begin{pmatrix}
                \lambda^{1}_{r_{1}+1} & & & & & & \\
                & \ddots & & & & & \\
                & & \lambda^{1}_{r}  & & & & \\
                & & & \ddots & & & \\
                & & & & \lambda^{n}_{r_{n}+1} & & \\
                & & & & & \ddots & \\
                & & & & & & \lambda^{n}_{r}
                        \end{pmatrix}\frac{\d \xi}{\xi}.
\end{equation}
Here, as in (1.11) of \cite{biqboa}, the constants 
$$
   \lambda^{j}_{k}=\frac{\mu^{j}_{k}-\beta^{j}_{k}}{2}
$$
are the eigenvalues of the residue at $p_j$ of the original Higgs field corresponding to 
$(E,\nabla,h)$, which are non-zero if and only if $\mu^{j}_{k}\neq 0$ because we assume 
that the parabolic connection is admissible and resonance-free. 

Recall from Lemma \ref{lem:ext} that the sections $\hat{\tau}^{\infty}_m(\xi,z)$ define a lattice 
$\widehat{F}$ near $\ti$, and that we set $\zeta=\xi^{-1}$ for a local coordinate of $\CPt$ near $\ti$. 
\begin{lmm}
With respect to the trivialization $\zeta\hat{\tau}^{\infty}_m(\xi,z)$ of the lattice 
$\widehat{F}\otimes \Hol_{\CPt}(\ti)$ the polar part of the transformed connection $\widehat{\nabla}$ reads 
$$
              - \begin{pmatrix}
                p_{1} & & & & & & \\
                & \ddots & & & & & \\
                & & p_{1}  & & & & \\
                & & & \ddots & & & \\
                & & & & p_{n} & & \\
                & & & & & \ddots & \\
                & & & & & & p_{n}
                        \end{pmatrix} \d \xi - 
               \begin{pmatrix}
                \mu^{1}_{r_{1}+1} & & & & & & \\
                & \ddots & & & & & \\
                & & \mu^{1}_{r}  & & & & \\
                & & & \ddots & & & \\
                & & & & \mu^{n}_{r_{n}+1} & & \\
                & & & & & \ddots & \\
                & & & & & & \mu^{n}_{r}
                        \end{pmatrix}\frac{\d \xi}{\xi}.
$$
\end{lmm}
\begin{proof}
Recall that the lattice $\widehat{F}$ had as generating local sections 
$\hat{\tau}^{\infty}_{k}(\xi)$ for all $p_j\in P$ and a diagonalizing basis $\tau^{j}_{k}$ of 
$\im(\res(\nabla,p_j))$. 
Now, (\ref{taukalapsigma}) states that the trivializations $\widehat{\mathcal{F}}$ and $\widehat{F}$ 
are linked (up to some multiplicative constants) as in formula (1.8) of \cite{biqboa}. 
By Theorem 4.32 of \cite{Sz}, the parabolic weight corresponding to the section 
$\hat{\sigma}^{\infty}_{k}(\xi)$ is $-1+\alpha^j_k$; in particular, all the weights 
corresponding to a trivialization of the lattice $\widehat{\mathcal{F}}$ lie in $[-1,0[$. 
Hence, the parabolic weights corresponding to $\widehat{\mathcal{F}}\otimes \Hol_{\CPt}(\ti)$ 
lie in $[0,1[$. Clearly, the trivializations $\zeta\hat{\sigma}^{\infty}_m(\xi,z)$ and 
$\zeta\hat{\tau}^{\infty}_m(\xi,z)$ of $\widehat{\mathcal{F}}\otimes \Hol_{\CPt}(\ti)$ 
and  $\widehat{F}\otimes \Hol_{\CPt}(\ti)$ 
respectively are then also linked as in formula (1.8) of \cite{biqboa}. 
Therefore, the polar parts of the transformed integrable connection $\widehat{\nabla}$ and of the 
transformed Higgs field $\widehat{\theta}$ with respect to the lattices 
$\widehat{\mathcal{F}}\otimes \Hol_{\CPt}(\ti)$ and $\widehat{F}\otimes \Hol_{\CPt}(\ti)$ are 
related as in formula (1.10) of loc. cit. On the other hand, in the Dolbeault theory 
the polar part does not change under tensoring a lattice by $\Hol_{\CPt}(\ti)$, so 
the polar part of the transformed Higgs field with respect to the lattice 
$\widehat{\mathcal{F}}\otimes \Hol_{\CPt}(\ti)$ is also equal to (\ref{stphDol}). 
This shows the lemma. 
\end{proof}

On the other hand, Lemma \ref{lem:parwght} implies that 
the parabolic weight corresponding to $\zeta\hat{\tau}^{\infty}_m(\xi,z)$ is $\beta^j_m$. 
This, however, means that the 
parabolic structure at $\ti$ is the direct sum of the parabolic structures at all 
the points $p_j\in P$ of the parabolic filtrations induced on $\im(A^j)$. 
As the filtration on $\coker(A^j)$ is by admissibility the canonical filtration, 
it follows that the parabolic structure at $\ti$ is the direct sum of the parabolic 
structures at all the points $p_j\in P$, ignoring the graded subspaces corresponding 
to the weight $0$.

\end{document}